\documentclass[a4paper,10pt]{article}
\usepackage[T1]{fontenc}
\usepackage[latin1]{inputenc}
\usepackage[dvips]{graphicx}
\usepackage[dvips]{color}
\usepackage{amssymb}
\usepackage{amsfonts}
\usepackage{amsmath}
\begin{document}


\title{Blow-up profiles of solutions for the exponential reaction-diffusion
equation.}
\author{A. Pulkkinen}
\date{}
\maketitle 

\begin{abstract}
We consider the blow-up of solutions for a semilinear reaction diffusion equation with exponential reaction term. It is know that certain solutions that can be continued beyond the blow-up time possess a nonconstant selfsimilar blow-up profile. Our aim is to find the final time blow-up profile for such solutions. The proof is based on general ideas using semigroup estimates. The same approach works also for the power nonlinearity.
\end{abstract}



\newtheorem{theo}{Theorem}
\newtheorem{pro}{Proposition}[section]
\newtheorem{lem}[pro]{Lemma}
\newtheorem{defi}[pro]{Definition}               
\newtheorem{coro}[pro]{Corollary}
\newtheorem{rema}[pro]{Remark}
\newtheorem{assu}[pro]{Assumption}

\newcommand{\ma}{\mathcal{A}}
\newcommand{\mb}{\mathcal{B}}
\newcommand{\mc}{\mathcal{C}}
\newcommand{\mf}{\mathcal{F}}
\newcommand{\mg}{\mathcal{G}}
\newcommand{\ml}{\mathcal{L}}
\newcommand{\mm}{\mathcal{M}}
\newcommand{\mn}{\mathcal{N}}
\newcommand{\mr}{\mathcal{R}}
\newcommand{\mt}{\mathcal{T}}
\newcommand{\mx}{\mathcal{X}}
\newcommand{\my}{\mathcal{Y}}
\newcommand{\mz}{\mathcal{Z}}                    
\newcommand{\mbc}{\mathbb{C}}
\newcommand{\mbn}{\mathbb{N}}
\newcommand{\mbz}{\mathbb{Z}}
\newcommand{\mbx}{\mathbb{X}}
\newcommand{\mbr}{\mathbb{R}}
\newcommand{\mbe}{\mathbb{E}}
\newcommand{\mbt}{\mathbb{T}}
\newcommand{\mbp}{\mathbb{P}}

\newcommand{\bdm}{\begin{displaymath}}
\newcommand{\edm}{\end{displaymath}}
\newcommand{\be}{\begin{equation}}
\newcommand{\ee}{\end{equation}}                   
\newcommand{\bmu}{\begin{multline*}}
\newcommand{\emu}{\end{multline*}}
\newcommand{\bea}{\begin{eqnarray*}}
\newcommand{\eea}{\end{eqnarray*}}

\newcommand{\comment}[1]{}
\newcommand{\eps}{\epsilon}
\newcommand{\veps}{\varepsilon}
\newcommand{\av}{\Arrowvert}
\newcommand{\ud}{\mathrm{d}}
\newcommand{\udm}{\mathrm{d} \mu}
\newcommand{\udmy}{\, \mathrm{d} \mu (y)}
\newcommand{\nelio}{\hspace{\stretch{1}}$\Box$}      
\newcommand{\wt}{\widetilde}
\newcommand{\wh}{\widehat}
\newcommand{\ol}{\overline}
\newcommand{\inti}{\int\displaylimits}
\newcommand{\lra}{\Leftrightarrow}
\newcommand{\llra}{\Longleftrightarrow}
\newcommand{\lto}{\longrightarrow}


\section{Introduction}

We consider the following problem
\begin{equation} \label{eq1}
\left\{
\begin{array}{lll}
u_t = \Delta u + f(u),& \quad x \in \Omega, &t>0, \\
u = 0,& \quad  x \in \partial \Omega, & t>0, \\
u(x,0) = u_0(x) \ge 0, & \quad x \in \Omega, &
\end{array}
\right.
\end{equation}
where $\Omega = B_R(0) = \{x \in \mbr^N \, : \, |x| < R\}$ and $N$ is supercritical, i.e., $N \in [3,9]$ and the initial condition $u_0$ is nonnegative and in $C^1(\ol{\Omega})$. We are mainly interested in the case $f(u) = e^u$, but some results work with more general nonlinearities. Before stating our results, see Theorems \ref{theorem_nondecrea}, \ref{theorem2} and \ref{theorem4} below, we give a brief introduction to the subject.

We are interested in solutions that blow up in finite time, which means that there exists $T\in (0,\infty)$ such that $\av u(\cdot,t) \av_{\infty} < \infty$ for $t < T$ and 
$$
\limsup_{t \to T}\av u(\cdot,t) \av_{\infty} = \infty.
$$
By standard theory of parabolic regularity, this implies that $u$ is a classical solution for every $t \in (0,T)$. Blow-up is said to be of type I if the blow-up rate is the same as that of the ordinary differential equation $u' = f(u)$. For $f(u) = e^u$ this means that 
$$
-C_1 \le \log(T-t) + \av u(\cdot,t) \av_{\infty} \le C_2,
$$
for some constants $C_1$ and $C_2$. If the blow-up is not of type I then it is said to be of type II. A point $x_0 \in \Omega$ is a blow-up point if there exists a sequence $\{(x_n,t_n)\}_n \subset \Omega \times (0,T)$ such that $(x_n,t_n) \to (x_0,T)$ and $u(x_n,t_n) \to \infty$ as $n \to \infty$.

A solution can exist beyond the blow-up time $t=T$ as a weak solution. To be more precise, we give the following definition.
\begin{defi}\label{d:L1}
By an $L^1$-solution of (\ref{eq1}) on $[0,\mt]$ we
mean a function $u\in C([0,\mt];$  $L^1(\Omega))$ such that $f(u)\in
L^1(Q_\mt)$, $Q_\mt:=\Omega\times(0,\mt)$ and the equality
$$\int_\Omega[u\Psi]^{t_2}_{t_1}\,dx-\int^{t_2}_{t_1} \int_\Omega u\Psi_t\,
dx\,dt=\int^{t_2}_{t_1} \int_\Omega (u\Delta\Psi+f(u)\Psi)\,dx\,dt$$
holds for any $0\le t_1<t_2\le \mt$ and $\Psi\in C^2(\bar Q_\mt)$,
$\Psi=0$ on $\partial\Omega\times [0,\mt]$.
\end{defi}
Blow-up is said to be complete if the solution can not be continued as an $L^1$-solution beyond the blow-up time.

In this paper we want to focus on solutions that blow up at $x = 0$ and have a nontrivial selfsimilar blow-up profile, i.e., the convergence (\ref{eq:ssconv1})-(\ref{stationaryAsympt}) below holds. The following theorem (proved in \cite{FP}) states that every radially symmetric solution that blows up and continues to exist as a weak solution has this property.
\begin{theo} \label{theo:FP}
Let $f(u) = e^u$, $N \in [3,9]$ and let $u$ be a radially symmetric and radially nonincreasing $L^1$-solution of (\ref{eq1}) on $[0,\mt]$ that blows up at $t =T < \mt$. Then 
\begin{equation} \label{eq:ssconv1}
\lim_{t \uparrow T} \big[ \log(T-t) + u(y\sqrt{T-t}, t)\big] = \varphi(y),
\end{equation}
uniformly for $y$ in compact sets of $\mbr^N$, where $\varphi$ satisfies
\begin{equation} \label{eq:stationary}
\left\{
\begin{array}{ll}
\Delta \varphi - \frac{y}{2} \nabla \varphi + e^\varphi -1 =0,& \qquad |y| >0, \\
\varphi(0) = \alpha, \; \nabla \varphi(0) = 0,
\end{array}
\right.
\end{equation}
and
\begin{equation} \label{stationaryAsympt}
\lim_{|y| \to \infty} \big( \varphi(y) +2\log |y| \big) = C_\alpha,
\end{equation}
for some $\alpha > 0$ and $C_{\alpha} \in \mbr$.
\end{theo}

If the above convergence (\ref{eq:ssconv1}) holds for some function $\varphi$, we will refer to $\varphi$ as the selfsimilar blow-up profile of $u$.

In the last section of this paper we will slightly improve the above Theorem \ref{theo:FP} by showing that the assumption on $u$ being radially nonicreasing is redundant. We will thus obtain the following.
\begin{theo} \label{theorem_nondecrea}
Let $u$ be a radially symmetric $L^1$-solution of (\ref{eq1}) with $f(u) = e^u$ on $[0,\mt]$ that blows up with type I rate at $(x,t) = (0,T)$, where $T < \mt$. Then the convergence (\ref{eq:ssconv1})-(\ref{stationaryAsympt}) holds.
\end{theo}

In this Theorem we assume that blow-up is of type I which holds if $N \in [3,9]$, $u$ is radially symmetric and the maximum of $u$ is attained at the origin, see \cite{FP}.

The existence of global $L^1$-solutions of (\ref{eq1}) with $f(u) = e^u$ that blow-up in finite time is proved  in \cite{LT} and \cite{FPo}. The previous two theorems then give the asymptotic behavior for such solutions as the blow-up time is approached. 

For subcritical dimensions $N \in [1,2]$ the only solution $\varphi$ of (\ref{eq:stationary})-(\ref{stationaryAsympt}) is $\varphi \equiv 0$ (see \cite{E}) and so the convergence (\ref{eq:ssconv1}) should always hold with $\varphi \equiv  0$. This is proved in \cite{BBE} under the extra assumption that $u$ is nondecreasing in time and radially decreasing. For general solutions the problem is how to obtain the blow-up rate.

For the power nonlinearity $f(u) = u|u|^{p-1}$ blow-up is always complete and $u$ has a constant selfsimilar blow-up profile, i.e., the covergence (\ref{eq:ssconv2}) holds for $\varphi$ equal to a constant, whenever $p < p_S$, where
$$
p_{S} = \left\{
\begin{array}{ll}
\infty,& \text{ if } N \le 2, \\
\frac{N+2}{N-2}, & \text{ if } N > 2,
\end{array} \right.,
$$
see \cite{GK}. 

For $f(u) = e^u$ and supercritical dimensions $N \in [3,9]$, however, there exists a sequence $\{\alpha_j\}$ of initial values tending to infinity such that the solutions $\varphi_j$ satisfy (\ref{eq:stationary})-(\ref{stationaryAsympt}) with $\alpha = \alpha_j$, see \cite{ET}. Similar result is also true for the power nonlinearity, \cite{BQ}.

The idea of the proof of Theorem \ref{theo:FP} is to assume that the convergence (\ref{eq:ssconv1}) holds with $\varphi \equiv 0$ and then prove that the final time blow-up profile of $u$ is given by
\be \label{loglog_profile}
u(x,T) + 2\log|x| + \log|\log|x|| \to C,
\ee
as $|x| \to 0$, which implies complete blow-up by the results in \cite{Va}. Convergence to a nontrivial $\varphi$ is then obtained by an energy argument. In \cite{B} the existence of solutions with the final time profile as in (\ref{loglog_profile}) is proved.

In \cite{HV1,Ve1,Ve2} the case of $f(u) = u^p$ is discussed and a variety of final time blow-up profiles is obtained by assuming constant selfsimilar blow-up profile. See also \cite{BB,FK} for other works in that direction. Later these methods were used in \cite{M} to prove Theorem \ref{theo:M} below, which corresponds to Theorem \ref{theo:FP}, but for the power type nonlinearity. The exponential nonlinearity in dimension one is discussed in \cite{HV2} and \cite{HV3}, and final time blow-up profiles are found, provided that the solution has a constant selfsimilar profile.

So there are many results concerning final time blow-up profiles of solutions provided that the selfsimilar blow-up profile is a constant one. The behavior of solutions as in Theorem \ref{theo:FP} at the blow-up moment is however not directly evident from the asymptotics (\ref{eq:ssconv1})-(\ref{stationaryAsympt}). Our main theorem of this paper is the following which in fact does give the final time blow-up profile for solutions satisfying (\ref{eq:ssconv1})-(\ref{stationaryAsympt}).

\begin{theo} \label{theorem2}
Assume that $u$ is a solution of (\ref{eq1}) with $f(u) = e^u$ that blows up with type I rate at $(x,t) = (0,T)$ for some $T<\infty$ and verifies (\ref{eq:ssconv1})-(\ref{stationaryAsympt}). Then the final time blow up profile of $u$ is given by
\be \label{eq:u_profile1}
|u(x,T) +2\log|x| - C_{\alpha}| \to 0, \quad \text{ as } |x| \to 0,
\ee
where $C_{\alpha}$ is the constant from (\ref{stationaryAsympt}).
\end{theo}

The existence of the limit $\lim_{t \to T} u(x,t)$ for $x \ne 0$ is a consequence of the parabolic estimates as will be seen in the proof of the above Theorem.

In this theorem we merely assume that $u$ is a continuous solution of (\ref{eq1}) that blows up with type I rate at $(x,t) = (0,T)$ and has a nontrivial selfsimilar blow-up profile, i.e., convergence as in (\ref{eq:ssconv1}) with (\ref{eq:stationary})-(\ref{stationaryAsympt}) holds. We do not need to assume that the solution is decreasing or even radially symmetric. It is of course a different matter whether there exist radially nonsymmetric solutions of (\ref{eq1}) verifying (\ref{eq:ssconv1})-(\ref{stationaryAsympt}) with radially symmetric $\varphi$. It is also not known if there exist any radially nonsymmetric selfsimilar solutions.

Even though the above Theorem \ref{theorem2} is stated with $f(u) = e^u$ our analysis works for a larger class of nonlinearities, including $f(u) = u^p$. For the algebraic nonlinearity Theorem \ref{theorem_nondecrea} corresponds to the following result, see \cite{M}.
\begin{theo} \label{theo:M}
Let $p>1$, $f(u) = u^p$ and $u$ be a radially symmetric $L^1$-solution of (\ref{eq1}) on $(0,\mt)$ that blows up with type I rate at $(x,t) = (0,T)$, where $T<\mt$. Then
\begin{equation}\label{eq:ssconv2}
\lim_{t \uparrow T} \big[(T-t)^{1/(p-1)} u(y\sqrt{T-t}, t)\big] = \varphi(y),
\end{equation}
uniformly for $y$ in compact sets of $\mbr^N$, where $\varphi$ satisfies
\begin{equation} \label{eq:stationary_pow}
\left\{
\begin{array}{ll}
\Delta \varphi -\frac{y}{2}\nabla \varphi - \frac{1}{p-1}\varphi + \varphi^p =0,& \qquad |y| >0, \\
\varphi(0) = \kappa + \alpha, \; \nabla \varphi(0) = 0,
\end{array}
\right.
\end{equation}
with
\be \label{kappa}
\kappa = \left(\frac{1}{p-1} \right)^{\frac{1}{p-1}}
\ee
and
\begin{equation} \label{stationaryAsympt2}
\lim_{|y| \to \infty} |y|^{2/(p-1)}\varphi(y)  = C_\alpha,
\end{equation}
for some $\alpha >0$ and $C_\alpha > 0$.
\end{theo}
Theorem \ref{theorem2} stated for $f(u) = u^p$ will then be the following, which is already known and proved in \cite{MM}. Our method, however, gives a new proof, which we do not present in this treatise, since it proceeds very much in the same way as the proof of Theorem \ref{theorem2} above.
\begin{theo} \label{theorem4}
Assume that $u$ is a solution of (\ref{eq1}) with $f(u) = u^p$ for some $p>1$ that blows up with type I rate at $(x,t) = (0,T)$ for some $T<\infty$ and verifies (\ref{eq:ssconv2})-(\ref{stationaryAsympt2}). Then the final time blow up profile of $u$ is given by
\be \label{eq:u_profile2}
\lim_{x\to 0} |x|^{2/(p-1)}u(x,T) = C_{\alpha},
\ee
where $C_{\alpha}$ is the same as in (\ref{stationaryAsympt2}).
\end{theo}

In a forthcoming paper \cite{P} we will show that if $u$ is a so-called minimal limit $L^1$-solution on $(0,\mt)$ that blows up at $t =T<\mt$ and if the assumptions of Theorem \ref{theorem_nondecrea} hold, then $u$ becomes regular immediately after the blow-up. Moreover, under some additional assumptions, the regularization is asymptotically selfsimilar, i.e., $u$ approaches a forward selfsimilar solution as $t \to T$ from above. This improves somewhat the results in \cite{FMP}.

The question about the behaviour of the final time blow-up profiles $u(x,T)$ near the blow-up point has been studied in many papers, but usually in the case where the selfsimilar blow up profile is the constant one, see \cite{HV1}, \cite{HV3}, \cite{Ve2}, \cite{M}. In these cases, the final time profile $u(x,T)$ is greater than those in (\ref{eq:u_profile1}) and (\ref{eq:u_profile2}) near the blow-up point and it holds that
$$
\lim_{x\to 0} |u(x,T) + 2\log|x|| = \infty,
$$
if $f(u) = e^u$, and
$$
\lim_{x\to 0} |x|^{2/(p-1)} u(x,T) = \infty,
$$
if $f(u) = u^p$. An example of this type of profiles for $f(u) = e^u$ is the one in (\ref{loglog_profile}).

Problems of this type have been studied by Matano and Merle in the paper \cite{MM} in more detail. The authors consider $f(u) = |u|^{p-1}u$, sign changing solutions, and $\Omega$ being either a ball or $\mbr^N$. Their result characterizes the size of the final time blow up profile in terms of the blow up rate and the behavior of the selfsimilar profile, but their technique does not seem to directly apply to the case with exponential nonlinearity and possibly increasing solutions.

The result of Matano and Merle gives that for $p > p_s$ one has
$$
\lim_{x \to 0} L^{-1}|x|^{2/(p-1)}u(x,T) =
\left\{
\begin{array}{lcl}
\infty \text{ or } -\infty & \Leftrightarrow & \text{ type I with } \varphi = \kappa \text{ or } -\kappa, \\
\text{finite but } \ne \pm 1, 0 & \Leftrightarrow & \text{ type I with nonconstant } \varphi, \\
1 \text{ or } -1 & \Leftrightarrow & \text{ type II}, \\
0 & \Leftrightarrow & \text{ no blow up at } x = 0,
\end{array}
\right.
$$
where $L^{p-1} = \frac{2}{p-1}\left(N-2 - \frac{2}{p-1}\right)$ and $\kappa$ is defined in (\ref{kappa}). Our Theorem \ref{theorem4} corresponds to the second equivalence.

Their techniques for obtaining the results are very different from ours. They obtain apriori bounds for the solutions and their derivatives by using some energy estimates and super solutions. They also work with radial solutions in order to be able to use parabolic estimates for one dimensional equations. These estimates then allow them to obtain the final time blow-up profile both in the case where $\varphi$ is regular and in the case where $\varphi$ is singular, and they prove immediate regularization with selfsimilar rate also for nonminimal $L^1$-solutions.

Our technique of proving Theorem \ref{theorem2} is based on the variation of constants formula and certain semigroup estimates. The assumptions here are not very strong and the ideas could be used also for different types of equations, but we cannot attack the situation where $\varphi$ is singular.

In paper \cite{P}, we use Theorem \ref{theorem4} to prove immediate regularization of solutions, but we can only consider the so-called minimal continuations.

Nonuniqueness of $L^1$-continuations of $u$ for $f(u) = u^p$ was proved in \cite{FM}.

The next section is devoted to a discussion on some properties of certain semigroups. We prove that the semigroup generated by the operator $\Lambda = \Delta - \frac{y}{2} \nabla + \Phi$ has specific regularization properties, similar to those of the semigroup generated by the Hermite operator $A=\Delta -\frac{y}{2} \nabla$, in case the function $\Phi=\Phi(y)$ decays to $0$ as $|y| \to \infty$.

In the third section we prove Theorem \ref{theorem2} by using the variation of constants formula, the semigroup estimates from Section 2, and some properties of the solution $u$ of (\ref{eq1}), specifically, the blow-up rate and the fact that
\be \label{eq:nablau}
|\nabla u(x,t) | \le \sqrt{2} e^{\max_x u(x,t)/2},
\ee
for every $x \in B_R(0)$ and $t \in (0,T)$, see \cite{FMc}.

In the last section we demonstrate that the results in \cite{FP} can be proved also without the assumption that $u$ is radially decreasing and thereby prove Theorem \ref{theorem_nondecrea}.

\section{Semigroup estimates}

To study the convergence (\ref{eq:ssconv1}) in more detail, we define the similarity variables $s = -\log(T-t)$ and $y = \frac{x}{\sqrt{T-t}}$ and let
\begin{equation} \label{eq:wdefinition_exp}
w(y,s) = \log(T-t) + u(x,t),
\end{equation}
for $|y| \le e^{s/2} R$ and $s \in [-\log(T),\infty)$. Here $u$ solves (\ref{eq1}) with $f(u) = e^u$. Then $w$ satisfies
\begin{equation} \label{eq:w}
w_s = \Delta w - \frac{y}{2} \nabla w + G(w), \quad \text{ for } |y| \le e^{s/2} R, \, s > -\log(T),
\end{equation}
with $w(y,s) = 0$ for $|y| = e^{s/2}R$, and
$$
G(w) = e^{w}-1.
$$
The convergence in (\ref{eq:ssconv1}) is equivalent to $w \to \varphi$ uniformly on compact sets as $s \to \infty$. Even though we take the equation with exponential nonlinearity to be the model problem, all the analysis goes through for a large class of nonlinearities.

Denoting $W = w- \varphi$ we notice that $W$ verifies
$$
W_s = \Delta W - \frac{y}{2}\nabla W + e^{\varphi} W + e^{\varphi}(e^W - 1 -W) = \Lambda W + e^{\varphi}(e^W-1-W),
$$
for $|y| \le R e^{s/2}$ and $s > -\log(T)$, where we have defined the operator
$$
\Lambda = \Delta -\frac{y}{2} \nabla + \Phi,
$$
with $\Phi = e^{\varphi}$.

The idea of the proof of Theorem \ref{theorem2} is simple. We prove that the convergence in (\ref{eq:ssconv1}) also holds for $|y| \le C (T-t)^{-1/2} = C e^{s/2}$. The claim is then achieved by using the asymptotics (\ref{stationaryAsympt}) of $\varphi$.

To obtain the convergence for $|y| \le C e^{s/2}$, we take a look at the shifted function $W(y + e^{s/2}\xi,s)$ as $s$ tends to infinity. Thus it is convenient to define the shifted and weighted $L^q$-norms as follows
\bdm
\mathcal{N}^q_r(\psi) = \sup_{|\xi| \le r} \Big( \int_{\mbr^N} |\psi(y)|^q e^{-(y-\xi)^2/4} \ud y \Big)^{1/q}, \quad \text{ for } r\ge 0,
\edm
and
\bdm
\ml^q_{\xi}(\psi) = \Big( \int_{\mbr^N} |\psi(y)|^q e^{-(y-\xi)^2/4} \ud y \Big)^{1/q}, \quad \text{ for } \xi \in \mbr^N.
\edm

In the following treatment we consider the semigroup generated by $\Lambda$ and assume only that $\Phi >0$ is bounded and verifies 
\be \label{eq:Phi_condition}
\Gamma = \max_{y \in \mbr^N} \Phi(y) < \infty \quad \text{and} \quad \Phi(y) \le \frac{C}{|y|^2},
\ee
for $y \in \mbr^N \setminus \{0\}$ and some constant $C>0$. We want to prove the necessary estimates that characterize the regularizing properties of the semigroup $\{e^{\Lambda t}\}_t$. Defining $A = \Delta - \frac{y}{2} \nabla$ to be the standard Hermite operator, we know that $A$ and $\Lambda $ are selfadjoint operators with domain $H^2_{\rho}(\mbr^N)$, which denotes the weighted Sobolev space with weight $\rho = \rho(y) = e^{-|y|^2/4}$. They generate strongly continuous semigroups in $L^2_{\rho}(\mbr^N)$, which we denote by $\{e^{At}\}_{t \ge 0}$ and $\{e^{\Lambda t}\}_{t \ge 0}$ respectively. We use the notation $\av \cdot \av_{L^2_{\rho}}$ for the norm in $L^2_{\rho}(\mbr^N)$.

We have the following formula for the action of the semigroup $e^{At}$ on functions $\psi \in L^2_{\rho}(\mbr^N)$,
\begin{equation} \label{eq:semigroup1}
e^{At} \psi(y) = \frac{1}{[4\pi(1-e^{-t})]^{N/2}} \int_{\mbr^N} \text{exp} \left( -\frac{(ye^{-t/2}-\lambda)^2}{4(1-e^{-t})} \right) \psi(\lambda) \ud \lambda, \quad \text{ for } t>0.
\end{equation}
Since the spectrum of $A$ consists of nonpositive real numbers we also know that there exists a constant $C>0$ such that
\begin{equation} \label{eq:semigroup2}
\av e^{At}\psi \av_{L^2_{\rho}} \le C \av \psi \av_{L^2_{\rho}},
\end{equation}
for every $\psi \in L^2_{\rho}(\mbr^N)$ and $t \ge 0$. Because we assumed (\ref{eq:Phi_condition}), it follows that
\begin{equation} \label{eq:semigroup3}
\av e^{\Lambda t}\psi \av_{L^2_{\rho}} \le C e^{\Gamma t} \av \psi \av_{L^2_{\rho}},
\end{equation}
for some constant $C>0$ and every $t \ge 0$.

The following regularizing property for the semigroup generated by $A$ can be found in \cite{Ve2}. The semigroup $\{e^{At}\}_t$ regularizes in the sense that it maps functions from $L^{\beta}_{\rho}(\mbr^N)$ to $L^q_{\rho}(\mbr^N)$ for any $q>\beta$ if $t$ is large enough. The following Proposition states this property in terms of the shifted norms.
\begin{pro} \label{pro:velaz}
Assume $1<q,\beta < \infty$ and $r, \wt{r} \ge 0$. Set $\beta' = \frac{\beta}{\beta-1}$. Then for any $t>0$ and any $\psi$ with  $\mn^q_{\wt{r}}(\psi)< \infty$, we have
\begin{multline} \label{L^q_r ineq1}
\mn^q_r(e^{At}\psi) \le \frac{1}{(4\pi(1-e^{-t}))^{N/2}} \Big( \frac{4\pi\beta(1-e^{-t})}{\beta'(\beta - 1 + e^{-t})}\Big)^{N/2\beta'}\\ \Big( \frac{4\pi(\beta - 1 + e^{-t})}{\beta - 1 - (q-1)e^{-t}} \Big)^{N/2q} \text{exp} \Big( \frac{e^{-t}(r - \wt{r}e^{t/2})_+^2}{4(\beta - 1 -(q-1)e^{-t})} \Big) \mn^{\beta}_{\wt{r}}(\psi),
\end{multline}
for $t \ge 0$ such that $\beta - 1 - (q-1)e^{-t} > 0$.
\end{pro}

We would like to obtain an analogous regularizing property for the semigroup generated by $\Lambda $.

This is done in different cases in Propositions \ref{pro:semigroup1} and \ref{pro:semigroup2} and Corollary \ref{coro:semigroup1} below. We will estimate the norm $\ml^q_{e^{s/2}\xi}(e^{\Lambda t}\psi)$ first when $t$ is strictly less than $s$ in Proposition \ref{pro:semigroup1}, then with $t=s$ in Proposition \ref{pro:semigroup2} and finally for $s$ large and $t$ close to $s$ in Corollary \ref{coro:semigroup1}.

Since $\frac{\ud}{\ud t} e^{\Lambda t} = \Lambda  e^{\Lambda t} = (A + \Phi )e^{\Lambda t}$, we know that
\begin{equation} \label{eq:variation}
e^{\Lambda t} \psi = e^{At} \psi + \int_0^t e^{A(t-\tau)} \Phi e^{\Lambda \tau}\psi \ud \tau.
\end{equation}
The next proposition restates Proposition \ref{pro:velaz} using the $\ml^q$-norms instead of $\mn^q$-norms.
\begin{pro} \label{pro:velaz2}
Assume that $1<q,\beta<\infty$ and $|\mu| \ge 0$. Then for any $\wh{t}$ such that $\frac{qe^{-\wh{t}}}{\beta - 1 + e^{-\wh{t}}} < 1$ there exists a constant $C$ such that
\bdm
\ml^q_{e^{t/2}\mu}(e^{At}\psi) \le C (1-e^{-t})^{-\eps} \ml^{\beta}_{\mu}(\psi),
\edm
for $\eps = \frac{N}{2\beta}$ and any $t>\wh{t}$. The constant $C$ is independent of $\mu$ and depends only on $\beta$, $q$, $N$ and $\wh{t}$.
\end{pro}
\emph{Proof.} From the proof of Proposition 2.1. in \cite{Ve2} one obtains that
$$
\ml^q_{\xi}(e^{At}\psi)^q \le C(\beta)(1-e^{-t})^{\frac{-Nq}{2\beta}} \int_{\mbr^N} \text{exp} \Big( -\frac{(y-\xi)^2}{4} + \frac{qe^{-t}(y-\mu e^{t/2})^2}{4(\beta - 1 + e^{-t})} \Big) \ud y \, \ml^{\beta}_{\mu} (\psi)^q.
$$
By taking $\xi = e^{t/2}\mu$ and by a change of variables we immediately get that
\begin{multline*}
\ml^q_{e^{t/2}\mu}(e^{At}\psi)^q \le  C(\beta)(1-e^{-t})^{\frac{-Nq}{2\beta}} \int_{\mbr^N} \text{exp} \Big( -\frac{y^2}{4} + \frac{qe^{-t}y^2}{4(\beta - 1 + e^{-t})} \Big) \ud y \, \ml^{\beta}_{\mu} (\psi)^q\\
= C(\beta)(1-e^{-t})^{\frac{-Nq}{2\beta}} \Big| \frac{qe^{-t}}{\beta - 1 + e^{-t}} - 1\Big|^{-N/2} \int_{\mbr^N} e^{-y^2/4} \ud y \, \ml^{\beta}_{\mu} (\psi)^q\\
= C(\beta,q,N,\wh{t})(1-e^{-t})^{\frac{-Nq}{2\beta}} \ml^{\beta}_{\mu} (\psi)^q
\end{multline*}
which gives the claim provided that $t>\wh{t}(q,\beta)$, because $\frac{\ud}{\ud t} \frac{qe^{-t}}{\beta-1+e^{-t}} < 0$. \nelio

Using this result we can prove regularizing properties also for the semigroup generated by the operator $\Lambda $. The first is the following.

\begin{pro} \label{pro:semigroup1}
Let $\Phi$ satisfy (\ref{eq:Phi_condition}). For every $\beta > N/2$, $0 < \delta < \beta-N/2$, and $r>0$ there exist constants $C = C(\delta,\beta,r,\Gamma)$, $\eps = \eps(\beta)>0$ and $M=M(\delta,\beta,r,\Gamma)>0$, such that
\be \label{eq:first_L_est}
\ml^p_{e^{s/2}\xi}(e^{\Lambda t}\psi) \le C (1-e^{-t})^{-\eps} \ml^{\beta}_{e^{(s-t)/2} \xi}(\psi),
\ee
which holds for every positive $t$ and $0 < t +M< s$, for every $p \in (1,\beta-\delta)$, every $|\xi| \ge  r$, and nonnegative $\psi\in L^{\beta}(\mbr^N)$.
\end{pro}

\emph{Proof.} Consider $\beta$, $\delta$, $r$, $\Phi$, $\xi$ and $\psi$ as in the claim. By the variation of constants formula,
\be \label{eq:estimate1}
 \ml^p_{e^{s/2}\xi}(e^{\Lambda t}\psi) \le \ml^p_{e^{s/2}\xi} (e^{At}\psi) + \int_0^{t} \ml^p_{e^{s/2}\xi}(e^{A(t-\tau)}\Phi e^{\Lambda \tau} \psi ) \ud \tau.
\ee
The previous Proposition \ref{pro:velaz2} implies that
\be \label{eq:ineq-1}
\sup_{p \in (1,\beta-\delta/3)} \ml^p_{e^{s/2}\xi}(e^{At}\psi) \le C_1(\delta,\beta)(1-e^{-t})^{-\eps} \ml^{\beta}_{e^{(s-t)/2}\xi}(\psi),
\ee
for some $C_1(\delta,\beta)>1$ and $\eps = \eps(\beta) = \frac{N}{2\beta} < 1$ for $t$ close to $0$.

Now we want to write estimates for the integral part in (\ref{eq:estimate1}). To that end, by using first the previous proposition and then Hölder's inequality, we get that there exists a constant $C_2(\delta,\beta)$ such that
\begin{multline} \label{eq:ineq0}
\ml^p_{e^{s/2}\xi}(e^{A(t-\tau)}\Phi e^{\Lambda \tau} \psi ) \le C_2(\delta,\beta) (1-e^{-(t-\tau)})^{-\eps'} \ml^{\beta'}_{e^{(s - t + \tau)/2} \xi}(\Phi e^{\Lambda \tau}\psi)
\\
\le C_2(\delta,\beta) (1-e^{-(t-\tau)})^{-\eps'} \ml^{\alpha}_{e^{(s-t+\tau)/2}\xi}(\Phi) \ml^{\beta''}_{e^{(s - t + \tau)/2} \xi}(e^{\Lambda \tau}\psi) ,
\end{multline}
for any $p \le \beta- \delta$ and $\beta' = \beta - 2\delta/3$, $\beta '' = \beta - \delta/3$, $\eps'=\frac{N}{2\beta'}$ and $\alpha = \beta' ( 1 - \frac{\beta'}{\beta''})^{-1}$. Because $|\Phi| \le \Gamma$ and verifies (\ref{eq:Phi_condition}), and since $|\xi| \ge r$, we can estimate, for any $\alpha>1$ and $t\ge 0$,
\begin{multline} \label{eq:phiestimate}
\ml^{\alpha}_{e^{t/2}\xi}(\Phi )^{\alpha}
\\
= \int_{|y| \le e^{t/2}|\xi|/2} \Phi(y)^{\alpha} e^{-(y-e^{t/2}\xi)^2/4} \ud y + \int_{|y| > e^{t/2}|\xi|/2} \ldots \ud y
\\
\le C_3 \Gamma^{\alpha} e^{Nt/2}|\xi|^N e^{-e^{t}|\xi|^2/16} + C_3 e^{-\alpha t} C_{\Phi}^{\alpha} |\xi|^{-2\alpha} \int_{|y|> e^{t/2} |\xi|/2} e^{-y^2/4} \ud y
\\
\le C_4(r,\Gamma,\alpha) e^{-\alpha t}.
\end{multline}

By the boundedness of $\Phi$, we can take $K>0$ such that, using Proposition \ref{pro:velaz2}, we have
\begin{multline} \label{eq:ineq1}
\ml^{\beta''}_{e^{s/2}\xi}(e^{\Lambda \tau}\psi) \le \ml^{\beta''}_{e^{s/2}\xi}(e^{(A+\Gamma)K/2} e^{\Lambda (\tau-K/2)}\psi)
\\
\le C_5(\Gamma,K)\ml^{\beta - \delta}_{e^{(s-K/2)/2}\xi}(e^{\Lambda (\tau-K/2)}\psi),
\end{multline}
for any $s > \tau \ge \frac{K}{2}$. By Proposition \ref{pro:velaz2}, we may also assume that the constant $C_5(\Gamma,K)>1$ is such that
\be \label{eq:ineq2}
\ml^{\beta''}_{e^{s/2}\xi}(e^{\Lambda \tau}\psi) \le \ml^{\beta''}_{e^{s/2}\xi}(e^{(A+\Gamma)\tau} \psi) \le C_5(\Gamma,K) (1-e^{-\tau})^{-\eps} \ml^{\beta}_{e^{(s-\tau)/2}\xi}(\psi),
\ee
for all $\tau < \frac{K}{2}$ and $s>\tau$. Using (\ref{eq:phiestimate}) and assuming that $M$ is large we obtain 
\begin{multline} \label{eq:ineq3}
\int_0^{t} \ml^{\alpha}_{e^{(s - t + \tau)/2}\xi}(\Phi) (1-e^{-(t-\tau)})^{-\eps'} (1-e^{-\tau})^{-\eps} \ud \tau
\\
\le \int_0^{\theta}\ldots \ud \tau + \int_{t-\theta}^{t} \ldots \ud \tau + C_6(\theta,r,\Gamma)\int_{\theta}^{t-\theta} e^{-(s- t + \tau)} \ud \tau (1-e^{-\theta})^{-\eps -\eps'}
\\
<  \frac{1}{2C_2(\delta,\beta)C_5(\Gamma,K)},
\end{multline}
for every $s \ge \frac M2$ and $t \in (0, s- \frac{M}{2})$. To show (\ref{eq:ineq3}), first take $\theta$ small and then $M$ large. Similarly, let $M$ be large enough so that also
\be \label{eq:ineq4}
\int_{K/2}^{t} \ml^{\alpha}_{e^{(s - t + \tau)/2}\xi}(\Phi) (1-e^{-(t-\tau)})^{-\eps'} (1-e^{-(\tau-K/2)})^{-\eps} \ud \tau < \frac{1}{8C_2(\delta,\beta) C_5(\Gamma,K)},
\ee
for every $s > t + M > K/2 + M$.

For $s > M$ let $t_B(s)$ be the supremum of such $t \in (0,s-M)$ for which
\bdm
\sup_{p \in (1,\beta-\delta)} \ml^p_{e^{s/2}\xi}(e^{\Lambda t}\psi) \le 2C_1(\delta,\beta)(1-e^{-t})^{-\eps} \ml^{\beta}_{e^{(s-t)/2}\xi}(\psi).
\edm
By inequality (\ref{eq:ineq-1}) and since the integral part in (\ref{eq:estimate1}) tends to zero as $t \to 0$, we know that $t_B(s)$ is positive. Let
$$
t(s) = \sup_{s'\in(M, s)} (s' - t_B(s')).
$$
We want to show that $t(s) \le M$ for every $s > M$, which implies $t_B(s) \ge s-M$ for every $s >M$. This will give the claim by the definition of $t_B(s)$.

Let us first show that $t_B(s) \ge \frac{K}{2}$. Assume, to obtain a contradiction, that $t_B(s) < \frac{K}{2}$ for some $s > M$. Then $t_B(s) < \frac{M}{2} < s - \frac{M}{2}$, since we may assume $M>K$. We may estimate the integral part of (\ref{eq:estimate1}) through
\begin{multline} \label{eq:intpart_est}
\int_0^{t_B(s)} \ml^p_{e^{s/2}\xi}(e^{A(t_B(s)-\tau)}\Phi e^{\Lambda \tau} \psi ) \ud \tau
\\
\le C_2 \int_0^{t_B(s)} (1-e^{-(t_B(s)-\tau)})^{-\eps'}  \ml^{\alpha}_{e^{(s - t_B(s) + \tau)/2} \xi} (\Phi) \ml^{\beta''}_{e^{(s-t_B(s)+\tau)/2}\xi}(e^{\Lambda \tau}\psi) \ud \tau
\\
< C_2 \int_0^{t_B(s)} (1-e^{-(t_B(s)-\tau)})^{-\eps'}  \ml^{\alpha}_{e^{(s - t_B(s) + \tau)/2} \xi} (\Phi) C_5 (1-e^{-\tau})^{-\eps} \ml^{\beta}_{e^{(s - t_B(s))/2}\xi}(\psi) \ud \tau
\\
< \frac{1}{2} \ml^{\beta}_{e^{(s - t_B(s))/2}\xi}(\psi),
\end{multline}
by inequalities (\ref{eq:ineq0}), (\ref{eq:ineq2}) and (\ref{eq:ineq3}). Here the constants $C_2$ and $C_5$ are as above, even though - for the sake of notation - their dependence on the parameters is not written out explicitly. Using this, together with (\ref{eq:estimate1}) and Proposition \ref{pro:velaz2}, we have
\begin{multline*}
\sup_{p \in (1,\beta-\delta)} \ml^p_{e^{s/2}\xi}(e^{\Lambda (t_B(s))}\psi)
\\
\le C_1(\delta,\beta)(1-e^{-t_B(s)})^{-\eps} \ml^{\beta}_{e^{(s - t_B(s))/2}\xi}(\psi) + \frac{1}{2} \ml^{\beta}_{e^{(s - t_B(s))/2}\xi}(\psi) \ud \tau
\\
\le \frac{3}{2} C_1(\delta,\beta)(1-e^{-t_B(s)})^{-\eps} \ml^{\beta}_{e^{(s - t_B(s))/2}\xi}(\psi),
\end{multline*}
by assuming $C_1(\delta,\beta)\ge 1$. This contradicts the definition of $t_B(s)$. Therefore $t_B(s) \ge \frac{K}{2}$ for every $s>M$ and so $t(s) \le s - K/2$.

Let us then show that $t(s) \le M$. We proceed again by contradiction and assume that $t(s) > M$ for some $s > M$, which implies $t(s) \in (M,s-K/2]$. Without loss of generality, we may assume that $t(s) = s- t_B(s)$, which gives $t_B(s) \in [K/2, s-M)$. For $\tau \in [\frac{K}{2},t_B(s)]$, we have that if $\wh{s} = t(s) + \tau - K/2$, then $\wh{s} - s = t(s) -s + \tau - K/2 \in [-t_B(s), -K/2]$ and so $\wh{s} < s$. Therefore $t(\wh{s}) \le t(s)$ and defining $\wh{\tau} = \tau - K/2$ we have $\wh{\tau} = \wh{s} - t(s) < \wh{s} - t(\wh{s}) \le t_B(\wh{s})$. Thus
\begin{multline} \label{eq:ineq5}
\ml^{\beta-\delta}_{e^{(s - t_B(s) + \tau- K/2)/2} \xi} (e^{\Lambda (\tau-K/2)}\psi) \le \sup_{p \in (1,\beta - \delta)} \ml^p_{e^{\wh{s}/2}\xi}(e^{\Lambda (\wh{\tau})})
\\
\le 2C_1(\delta,\beta)(1-e^{-(\tau - K/2)})^{-\eps}\ml^{\beta}_{e^{t(s)/2}\xi}(\psi),
\end{multline}
by the definition of $t_B(\wh{s})$.

Precisely as in (\ref{eq:intpart_est}), we obtain
\bdm
\int_0^{K/2} \ml^p_{e^{s/2}\xi}(e^{A(t_B(s)-\tau)}\Phi e^{\Lambda \tau} \psi ) \ud \tau  < \frac{1}{2} \ml^{\beta}_{e^{(s - t_B(s))/2}\xi}(\psi)
\edm
and, by inequalities (\ref{eq:ineq0}), (\ref{eq:ineq1}), (\ref{eq:ineq4}) and  (\ref{eq:ineq5}), we can estimate
\begin{multline*}
\int_{K/2}^{t_B(s)} \ml^p_{e^{s/2}\xi}(e^{A(t_B(s)-\tau)}\Phi e^{\Lambda \tau} \psi ) \ud \tau
\\
\le C_2\int_{K/2}^{t_B(s)} \ml^{\alpha}_{e^{(s - t_B(s) + \tau)/2} \xi} (\Phi )(1-e^{-(t_B(s)-\tau)})^{-\eps'} \ml^{\beta''}_{e^{(s - t_B(s) + \tau)/2} \xi} (e^{\Lambda \tau}\psi) \ud \tau
\\
\le C_2 C_5\int_{K/2}^{t_B(s)} \ml^{\alpha}_{e^{(s- t_B(s) + \tau)/2} \xi} (\Phi )(1-e^{-(t_B(s)-\tau)})^{-\eps'} \phantom{jotainjotainjotianjo}
\\
\phantom{jotaintajoatjotaojoajt}\cdot \ml^{\beta-\delta}_{e^{(s- t_B(s) + \tau - K/2)/2} \xi} (e^{\Lambda (\tau-K/2)}\psi) \ud \tau
\\
\le 2C_1 C_2 C_5 \int_{K/2}^{t_B(s)} \ml^{\alpha}_{e^{(s - t_B(s) + \tau)/2} \xi} (\Phi )(1-e^{-(t_B(s)-\tau)})^{-\eps'} (1-e^{-(\tau -K/2)})^{-\eps} \ud \tau
\\
\phantom{jotaintjojojojojojojottttttttttttttttttttttttt} \cdot \ml^{\beta}_{e^{(s - t_B(s))/2}}(\psi)
\\
\le \frac{1}{4}C_1 \ml^{\beta}_{e^{(s - t_B(s))/2}}(\psi) <\frac{1}{4}C_1(1-e^{-t_B(s)})^{-\eps} \ml^{\beta}_{e^{(s-  t_B(s))/2}\xi}(\psi),
\end{multline*}
where the dependence of $C_1$, $C_2$ and $C_5$ on the parameters is not explicitly written out.

This implies that
\bdm
\sup_{p \in (1,\beta-\delta)} \ml^p_{e^{s/2}\xi}(e^{\Lambda (t_B(s))}\psi) \le (1 + \frac{1}{2} + \frac{1}{4} )C_1(\delta,\beta)(1-e^{-t_B(s)})^{-\eps} \ml^{\beta}_{e^{(s - t_B(s))/2}\xi}(\psi),
\edm
which contradicts the definition of $t_B(s)$. Therefore the only possibility is that $t(s) \le M$ and so the claim is proved.  \nelio

In the previous proposition we assumed that $t < s - M$ for some large $M$. The next Proposition deals with the case $t = s > M$ for some large $M$.

\begin{pro} \label{pro:semigroup2}
Let $\Phi$ satisfy (\ref{eq:Phi_condition}). For every $r\in (0,1)$ there exist constants $M = M(r,\Gamma)$ and $C = C(r,\Gamma)$ such that 
\bdm
\ml^2_{e^{s/2}\xi}(e^{\Lambda s}\psi) \le C(r,\Gamma,N) \av \psi \av_{L^2_{\rho}}, 
\edm
for every $s > M$, $|\xi| \in (r,1/r)$, and every nonnegative $\psi \in L^2_{\rho}(\mbr^N)$. 
\end{pro}

\emph{Proof.} Fix $r>0$ and take $\theta = \log(2N)$. Assume that $s\ge M + \theta $, where $M>\theta$ will be defined later and let $|\xi| \in (r,\frac{1}{r})$.

Let us first note that for any $1<p<q<\infty$ and $\frac{p}{q} + \frac{1}{\gamma} = 1$, we have
\begin{multline} \label{pq_norm_est}
\ml^p_{\xi}(\psi)^p = \int_{\mbr^N} |\psi|^p e^{-p|y|^2/4q} e^{(1-\frac{1}{\gamma})|y|^2/4} e^{-(y-\xi)^2/4} \ud y
\\
\le \left( \int_{\mbr^N} |\psi|^q e^{-|y|^2/4} \ud y \right)^{p/q} \left( \int_{\mbr^N} e^{(\gamma-1)|y|^2/4} e^{-\gamma(y-\xi)^2/4} \ud y \right)^{1/\gamma} 
\\
= \av \psi \av_{L^q_{\rho}}^p \left( \int_{\mbr^N} e^{-\frac{y}{2}(\frac{y}{2} - \gamma \xi) - \frac{\gamma}{4} |\xi|^2} \ud y \right)^{1/\gamma} \le \av \psi \av_{L^q_{\rho}}^p \left( \int_{\mbr^N} e^{-(\frac{y}{2}+ \frac{\gamma}{2}\xi)(\frac{y}{2} - \frac{\gamma}{2} \xi) - \frac{\gamma}{4} |\xi|^2}  \ud y\right)^{1/\gamma}
\\
= \av \psi \av_{L^q_{\rho}}^p \left( \int_{\mbr^N} e^{-|y|^2/4} e^{(\gamma^2-\gamma)|\xi|^2/4} \ud y  \right)^{1/\gamma} = C(p,q,r)  \av \psi \av_{L^q_{\rho}}^p
\end{multline}

If $t \le  M$, then $s-t \ge \theta$ and so by using first Proposition \ref{pro:velaz2}, then Hölder's inequality, then the fact that $\ml^{21/2}_{e^{t/2}\xi}(\Phi ) < C(r,\Gamma) e^{-t}$ by (\ref{eq:phiestimate}), and finally the above inequality (\ref{pq_norm_est}), we have
\begin{multline} \label{eq:pro2.4.1}
\ml^2_{e^{s/2}\xi}(e^{A(s-t)} \Phi  e^{\Lambda t} \psi) \le C_1(\theta) \ml^{3/2}_{e^{t/2} \xi}(\Phi e^{\Lambda t} \psi) \le C_1(\theta)  \ml^{21/2}_{e^{t/2}\xi}(\Phi) \ml^{7/4}_{e^{t/2}\xi}(e^{\Lambda t}\psi)
\\
\le C_2(\theta,r,\Gamma) e^{-t} e^{\Gamma M} \mn^{7/4}_{e^{M/2}\xi} (e^{A t}\psi) \le C_3(\theta,r,\Gamma,M) e^{-t} \av e^{At} \psi \av_{L^2_{\rho}}
\\
\le C_3(\theta,r,\Gamma,M) e^{-t} \av \psi \av_{L^2_{\rho}}.
\end{multline}

If  $0<t<s-\theta$, then because of our choice of $\theta$, we can use Proposition \ref{pro:velaz2} with $q = 2$ and $\beta = \frac{3}{2}$ for the first inequality below and Hölders inequality for the second, to obtain
\begin{multline} \label{eq:pro2.4.2}
\ml^2_{e^{s/2}\xi}(e^{A(s-t)} \Phi  e^{\Lambda t} \psi) \le C_1(\theta) \ml^{3/2}_{e^{t/2}\xi}(\Phi e^{\Lambda t}\psi)
\\
\le C_1(\theta)\ml^{12/2}_{e^{t/2}\xi}(\Phi) \ml^{2}_{e^{t/2}\xi}(e^{\Lambda t}\psi) \le C_4(\theta,r,\Gamma) e^{-t} \ml^{2}_{e^{t/2}\xi}(e^{\Lambda t}\psi),
\end{multline}
where the last inequality is again due to (\ref{eq:phiestimate}).

If $M \le s-\theta \le t \le s$, then using first Proposition \ref{pro:velaz2} with exponents $q = 2$ and $\beta =N$, then Hölder's inequality, and finally Proposition \ref{pro:velaz2} with exponents $q = 2N$ and $\beta = 2$, we obtain
\begin{multline} \label{eq:pro2.4.3}
\ml^2_{e^{s/2}\xi} (e^{A(s-t)} \Phi  e^{\Lambda t} \psi) \le C_5(1-e^{-(s-t)})^{-1/2} \ml^{N}_{e^{t/2}\xi}( \Phi e^{\Lambda t} \psi)
\\
\le C_5 \ml^{2N}_{e^{t/2}\xi}(\Phi )(1-e^{-(s-t)})^{-1/2} \ml^{2N}_{e^{t/2}\xi}(e^{(A+\Gamma)\theta} e^{\Lambda (t-\theta)}\psi)
\\
\le C_6(\theta,r,\Gamma)e^{-t} (1-e^{-(s-t)})^{-1/2} \ml^{2}_{e^{(t-\theta)/2}\xi}(e^{\Lambda (t-\theta)}\psi).
\end{multline}

Define $W(s) = \ml^2_{e^{s/2}\xi}(e^{\Lambda s}\psi)$ and use the variation of constants formula (\ref{eq:variation}) together with Proposition \ref{pro:velaz} and the above estimates (\ref{eq:pro2.4.1}) - (\ref{eq:pro2.4.3}) to observe that
\begin{multline} \label{W(s)_est}
W(s) \le \ml^2_{e^{s/2}\xi}(e^{As} \psi) + \int_0^s \ml^2_{e^{s/2}\xi} (e^{A(s-t)}\Phi  e^{\Lambda t}\psi) \ud t \\
\le  C_7 \av \psi \av + \int_0^{M} C_3 e^{-t} \av \psi \av_{L^2_{\rho}} \ud t \\
+ \int_{M}^{s-\theta} C_4 e^{-t} W(t) \ud t \\
+ \int_{s-\theta}^s C_6 e^{-t} (1-e^{-(s-t)})^{-1/2} W(t-\theta) \ud t 
\end{multline}
for every $s \ge M + \theta$ and for some $C_7>0$ arising from the estimate (\ref{L^q_r ineq1}).

By Proposition \ref{pro:velaz}, we have, for every $s \in [M,M+2\theta]$, the estimate
\be \label{M,M2t_est}
W(s) \le \ml^2_{e^{s/2}\xi}(e^{(A+\Gamma)s}\psi) \le e^{\Gamma(M+2\,\theta)} \mn^2_{e^{s/2}\xi} (e^{A s} \psi) \le C_8 e^{\Gamma(M+2\theta)} \av \psi \av_{L^2_{\rho}}.
\ee
Take $M = M(\theta,r,\Gamma)$ large enough such that
\be \label{M_def}
 C_4(\theta,r,\Gamma) e^{-M} + C_6(\theta,r,\Gamma) e^{-M} \int_0^{\theta} (1-e^{-t}) \ud t < \frac{1}{3}.
\ee
Then let $\wt{C}(\theta,r,\Gamma) = C_8 e^{\Gamma(M+\theta)}$ and take $K = K(\theta,r,\Gamma)>1$ large enough such that
\be \label{K_def}
\frac{C_7}{K\wt{C}(\theta,r,\Gamma)} + \frac{ C_3(\theta,r,\Gamma,M)}{K\wt{C}(\theta,r,\Gamma)} < \frac{1}{3}.
\ee
Let $\wt{s}$ be the supremum of such $s' > M$ for which $W(s) \le K \wt{C}(\theta,r,\Gamma) \av \psi \av_{L^2_{\rho}}$ for every $s \in [M,s']$. By (\ref{M,M2t_est}) we know that $\wt{s} \ge M + 2\, \theta$. Assuming that $\wt{s} < \infty$, by using (\ref{W(s)_est}), (\ref{M_def}) and (\ref{K_def}), we have for every $s \in [M+2\theta,\wt{s}]$ that
\begin{multline*}
W(s) \le C_7 \av \psi \av_{L^2_{\rho}}  + C_3 \av \psi \av_{L^2_{\rho}} 
\\
+ C_4 (e^{-M} - e^{-(s-\theta)}) K \wt{C}(\theta,r,\Gamma) \av \psi \av_{L^2_{\rho}} 
\\
+ C_6 e^{-M} K \wt{C}(\theta,r,\Gamma)  \int_0^{\theta} (1-e^{-t}) \ud t \av \psi \av_{L^2_{\rho}} 
\\
\le K \wt{C}(\theta,r,\Gamma) \av \psi \av_{L^2_{\rho}} \Big( \frac{C_7}{K\wt{C}(\theta,r,\Gamma)} +\frac{C_3}{K\wt{C}(\theta,r,\Gamma)} +  C_4e^{-M}
\\
+  C_6 e^{-M}\int_0^{\theta} (1-e^{-t}) \ud t \Big) \le \frac{2}{3} K \wt{C}(\theta)  \av \psi \av_{L^2_{\rho}},
\end{multline*}
which contradicts the definition of $\wt{s}$ and so $\wt{s} = \infty$. We have therefore obtained that 
$$
W(s) = \ml^2_{e^{s/2}\xi}(e^{\Lambda s}\psi) \le K(\theta,r,\Gamma) \wt{C}(\theta,r,\Gamma) \av \psi \av_{L^2_{\rho}},
$$
for every $s > M(\theta,r,\Gamma)$. This gives the claim. \nelio

The following corollary is almost a restatement of the previous proposition, but instead of considering $\ml^2_{e^{s/2}\xi}(e^{\Lambda t} \psi)$ for $s = t$, we allow $s-t$ to be positive and bounded.

\begin{coro}\label{coro:semigroup1}
Let $\Phi$ satisfy (\ref{eq:Phi_condition}). For every $r\in(0,1)$, and $M>0$ there exist constants $K = K(M,r,\Gamma)$ and $C = C(M,r,\Gamma)$ such that
\bdm
\ml^2_{e^{(s-s_0)/2}\xi}(e^{\Lambda (s-t)} \psi) \le C(M,r,\Gamma) \av \psi \av_{L^2_{\rho}},
\edm
for every $t \in [s_0, s_0 + M]$, $s > s_0 + K(M,r,\Gamma)$, $s_0 \ge 0$, every $|\xi| \ge (r,\frac 1r )$, and every nonnegative $\psi \in L^2_{\rho}(\mbr^N)$.
\end{coro} 

\emph{Proof.} By the previous Proposition \ref{pro:semigroup2} there exist constants $M'(r,\Gamma)$ and $C(M,r,\Gamma)$ such that
\bdm
\ml^2_{e^{\tau/2} \eta}(e^{\Lambda  \tau}\psi) \le C(M,r,\Gamma) \av \psi \av_{L^2_{\rho}},
\edm
for $\tau > M'(r,\Gamma)$ and $|\eta| \in (r, e^M/r )$.

Assume that $s,t$ and $\xi$ are as in the claim, and let $\tau = s-t$ and $\eta = e^{t-s_0} \xi$. Then $\tau > K(M,r,\Gamma) - M > M'(r,\Gamma)$ by choosing $K$ large enough, and $|\eta| \in (r,e^M/r)$, which implies that
\bdm
\ml^2_{e^{(s-s_0)/2}\xi}(e^{\Lambda (s-t)}\psi) = L^2_{e^{\tau/2} \eta}(e^{\Lambda  \tau}\psi) \le C(M,r,\Gamma) \av \psi \av_{L^2_{\rho}}, 
\edm
and so the claim is proved. \nelio

We have now the desired results that describe the regularizing properties of the semigroup generated by $\Lambda$. These results will be used in the next section to prove Proposition \ref{pro:final}.

\section{The final time blow-up profile}

In this section we prove Theorem \ref{theorem2}.

Let $u$ and $\varphi$ be as in Theorem \ref{theorem2} and define the usual similarity variables through $s = -\log(T-t)$ and $y = \frac{x}{\sqrt{T-t}}$ and let 
$$
W(y,s) = \log(T-t) + u(x,t) - \varphi(y).
$$
Our assumptions in Theorem \ref{theorem2} imply that $W(y,s) \to 0$ uniformly for $y$ in compact sets as $s \to \infty$.

Generally speaking, to prove Theorem \ref{theorem2}, we want to show that the $\ml^2_{e^{(s-s_0)/2}\xi}$-norm of $W(\cdot,s)$ can be estimated by the $L^2_{\rho}$-norm of $W(\cdot,s_0)$. This is done in Proposition \ref{pro:final} below. In the proof of that Proposition we utilize the regularizing properties of the semigroup $\{e^{\Lambda t}\}_t$, obtained in the previous section. Using this results and the $L^{2}$ - $L^{\infty}$ regularization of the semigroup generated by $A$, one has that $|W(e^{(s-s_0)/2}\xi,s)|$ tends to zero as $s_0$ tends to infinity. By the definition of $W$ and by the asymptotics (\ref{stationaryAsympt}) of $\varphi$ we obtain the claim at the very end of this section.

Before stating and proving Proposition \ref{pro:final}, we have to consider the properties of the function $W$ in more detail. We also need some auxiliary results. In Proposition \ref{pro:wpro1} we demonstrate how to move from $\ml^2_{\xi}$ -norm to $\mn^2_{|\xi|}$ -norm, in Proposition \ref{pro:boundedgrad2} we consider the $L^2$ - $L^{\infty}$ regularization of $e^{At}$ and in Proposition \ref{pro:hpro} we estimate the norm of the nonlinearity appearing in the equation for $W$.

Since the function $W$ is defined on some $s$ dependent subset of $\mbr^N$, we need to extend it to $\mbr^N$. Because the blow-up set is a compact set of $B(R)$, we can take $R_1 \in (0,R)$ such that $u(x,t)$ is bounded for $(x,t) \in B(R)\setminus B(R_1) \times (0,T)$. Then let $\zeta$ be a smooth function equal to $1$ for $|x| \le R_1$ and equal to $0$ for $|x|>R$ and let
\bdm
\wt{W}(y,s) = \zeta(e^{-s/2}y) W(y,s).
\edm

Now $\wt{W}$ is defined in the whole space $\mbr^N$ and it satisfies the equation
\bdm
\wt{W}_s = \Lambda  \wt{W} + \wt{h},
\edm
where
\begin{equation*}
\wt{h}= -e^{-s} \Delta\zeta W - 2e^{-s/2} \nabla \zeta \cdot \nabla W \\
+ \frac{e^{-s/2}y}{2}\cdot \nabla \zeta W + \zeta e^{\varphi}(e^{W} - 1) - e^{\varphi} \wt{W},
\end{equation*}
for $|y| \le R e^{s/2}$ and $\wt{h} \equiv 0$, for $|y| > R e^{s/2}$. Here we use the notation $\nabla \zeta$ for $\nabla \zeta(e^{s/2}y)$ and the same applies to the Laplacian of $\zeta$.

Since the blow-up is assumed to be of type I, we obtain, by using (\ref{eq:nablau}),
$$
|\nabla \wt{W}(y,s)| = |\nabla W(y,s)| = |\sqrt{T-t}\nabla u(x,t) - \nabla \varphi(y)|\le C,
$$
for $|y| \le R_1 e^{s/2}$. For such $y$ one also has $\wt{h} = e^{\varphi} (e^{\wt{W}} - 1 - \wt{W})$. Moreover, type I blow-up implies that $W$ is bounded from above and so the estimates $|\wt{h}| \le C |\wt{W}|$ and $|\wt{h}| \le C |\wt{W}|^2$ are valid for some constant $C>0$.

For $|y| = |x| e^{s/2}$ and $|x| \in (R_1,R)$, one has
\begin{multline*}
|\nabla \wt{W}(y,s)| 
\\
\le e^{-s/2} |\nabla \zeta(e^{-s/2}y)|\big(u(x,t) + |s + \varphi(y)|\big) + |\zeta(e^{-s/2}y)||\sqrt{T-t}\nabla u(x,t) -\nabla \varphi(y)| \le C,
\end{multline*}
by using (\ref{eq:nablau}) and the asymptotic behavior (\ref{stationaryAsympt}) of $\varphi$. Similarly we have
$$
|\wt{W}(y,s)| \le |W(y,s)| \le u(x,t) + |\varphi(e^{s/2}x)+s| \le C.
$$
Therefore, we obtain the estimate
$$
|\wt{h}| \le A_0,
$$
for some finite constant $A_0>0$.

For $|y| > Re^{s/2}$ we have that $\wt{W} = 0$ and $\wt{h} = 0$.

Defining
$$
Z(y,s) = |\wt{W}(y,s)| = |\zeta(e^{-s/2}y)( -s + u(e^{-s/2}y,T-e^{-s}) - \varphi(y)|,
$$
we have obtained that $Z$ satisfies
\be \label{Z_equation}
Z_s \le \Lambda Z + h + A_0 \chi,
\ee
where $\chi(y,s) = \chi_{\{|y| > R_1 e^{s/2}\}}(y)$ is the characteristic function of the set $\{y \in \mbr^N \, : \, |y| > R_1 e^{s/2} \}$. The function $h$ verifies
\be \label{esti1}
h \le A_1 Z, \quad \text{for } y \in \mbr^N,
\ee
and, for some $\vartheta(y,s) \in (0,\wt{W}(y,s))$,
\be \label{esti2}
h = \frac{1}{2} e^{\varphi+ \vartheta} Z^2 \le A_2 Z^2, \quad \text{for } |y| \le R_1 e^{s/2},
\ee
for some constants $A_1$ and $A_2$. Since $\varphi$ is assumed to be as in Theorem \ref{theo:FP}, we have that
$$
\max_{y \in \mbr^N} \varphi(y) = \varphi(0) = \alpha.
$$
Therefore $Z$ satisfies also the inequality
\be \label{Z_equation_2}
Z_s \le (A+\Gamma+A_1)Z + A_0 \chi,
\ee
where $\Gamma = e^{\alpha}$. Since $\wt{W}(y,s) \to 0$ on compact sets as $s \to \infty$ and since $|\nabla \wt{W}| \le C$, one has
\be \label{esti3}
\av Z(\cdot,s) \av_{L^2_{\rho}} \le A_3, \quad \text{for }s\ge 0 \text{ and} \quad \av Z(\cdot,s) \av_{L^2_{\rho}} \to 0,
\ee
as $s \to \infty$. Moreover,
\be \label{nablaZ}
|\nabla Z(\cdot,s)| \le A_4, \quad \text{for }(y,s) \in \mbr^N \times (0,\infty).
\ee

In what follows, we will consider the parameters $\{A_i\}_{i=0}^4$, $R_1$ and $\Gamma$ as given. The constants below may depend on these parameters, but we will not state it explicitly.

Let us now derive some estimates for the $\ml^q_{\xi}$ -norm of $Z$. To that end, let $t \in (s-\theta,s)$, $|\lambda| > R_1 e^{t/2}$ and $|y| < R_1 e^{s/2}/2$, which gives
$$
\frac{|ye^{-(s-t)/2} - \lambda|}{\sqrt{1-e^{-s+t}}} > \frac{e^{t/2} R_1}{2\sqrt{1-e^{-\theta}}} > \frac{e^{t/2}R_1}{2},
$$
and so, by using the representation formula (\ref{eq:semigroup1}) for the semigroup, we have that, for such $t$ and $y$, it holds
\begin{multline} \label{ineq:Achi_1}
e^{A(s-t)}\chi(y,t)
\\
= \frac{1}{(4\pi(1-e^{-s+t}))^{N/2}} \int_{\mbr^N} \exp \left( - \frac{(ye^{-(s-t)/2} - \lambda)^2}{4(1-e^{-s+t})} \right)\chi(\lambda,t) \ud \lambda
\\
= \frac{1}{(4\pi(1-e^{-s+t}))^{N/2}} \int_{|\lambda| > R_1e^{t/2}} \exp \left( - \frac{(ye^{-(s-t)/2} - \lambda)^2}{4(1-e^{-s+t})} \right) \ud \lambda
\\
\le \frac{1}{(4 \pi)^{N/2}} \int_{|\lambda| > \frac{e^{t/2} R_1}{2}} e^{-|\lambda|^2/4} \ud \lambda \le C_1 e^{-t}.
\end{multline}
Above, to be more precise, we could have written $e^{A(s-t)}\chi(y,t)$ as $[e^{A(s-t)}\chi(\cdot,t)](y)$, but we obey the former option in what follows. If $t \in (s-\theta,s)$ and $|y| > R_1 e^{s/2}/2$, then
\begin{multline*}
e^{A(s-t)}\chi(y,t) \le \frac{1}{(4\pi(1-e^{-s+t}))^{N/2}} \int_{\mbr^N} \exp \left( - \frac{(ye^{-(s-t)/2} - \lambda)^2}{4(1-e^{-s+t})} \right) \ud \lambda 
\\
= \frac{1}{(4\pi)^{N/2}}\int_{\mbr^N} e^{-|\lambda|^2/4} \ud \lambda = C_2.
\end{multline*}
Therefore, for any $\sigma>0$, $|\xi| \le \max \{\sigma,\frac{e^{s/2}R_1}{4}\}$, and $t \in (s-\theta,s)$, we have
\begin{multline*}
\ml^q_{\xi}(e^{A(s-t)} \chi(\cdot,t))^q 
\\
\le C_1^q e^{-qt} \int_{|y| < R_1 e^{s/2}/2} e^{-(y-\xi)^2/4} \ud y + C_2^q \int_{|y| > R_1 e^{s/2}/2} e^{-(y-\xi)^2/4} \ud y 
\\
\le C_3(\sigma,q)^q e^{-qt}.
\end{multline*}
Next we use the variation of constants formula and (\ref{Z_equation_2}) to obtain
\begin{multline} \label{ml_Z_esti_0}
\ml^q_{\xi}(Z(\cdot,s)) \le \ml^q_{\xi}(e^{(A+\Gamma+A_1)\theta}Z(\cdot,s-\theta)) + \int_{s-\theta}^s \ml^q_{\xi}(e^{(A+\Gamma+A_1)(s-t)}A_0 \chi(\cdot,t) ) \ud t
\\
= \ml^q_{\xi}(e^{(A+\Gamma+A_1)\theta}Z(\cdot,s-\theta)) + A_0 C_3(\sigma,q) \int_{s-\theta}^s e^{(\Gamma + A_1)(s-t)}e^{-t}\ud t
\\
= \ml^q_{\xi}(e^{(A+\Gamma+A_1)\theta}Z(\cdot,s-\theta)) + \frac{A_0 C_3(\sigma,q) e^{-s}}{\Gamma + A_1 + 1 }(e^{(\Gamma + A_1 + 1 )\theta} - 1),
\end{multline}
for every $|\xi| \le \max\{\sigma,\frac{e^{s/2}R_1}{4}\}$. By Proposition \ref{pro:velaz2}, we then get
\be \label{ml_Z_esti}
\ml^q_{e^{\theta/2}\mu}(Z(\cdot,s)) \le C_4(\sigma,q,\beta,\theta) \left( \ml^{\beta}_{\mu}(Z(\cdot,s-\theta))+e^{-s} \right),
\ee
for every $|e^{\theta/2}\mu| \le \max\{\sigma,\frac{e^{s/2}R_1}{4}\}$, $s > \theta$, and $\theta$ such that $\frac{q e^{-\theta}}{\beta-1+e^{-\theta}} < 1$. By Proposition \ref{pro:velaz} it also holds that
\begin{multline} \label{mn_Z_esti}
\mn^2_{\sigma}(Z(\cdot,s)) \le e^{(\Gamma+A_1)\theta} \mn^2_{\sigma}(e^{A \theta} Z(\cdot,s-\theta)) + C_5'(\theta,\sigma)e^{-s}
\\
\le C_5(\theta,\sigma) \left( \av Z(\cdot,s-\theta) \av_{L^2_{\rho}} + e^{-s} \right),
\end{multline}
for any $\theta, \sigma > 0$.

We want to consider the solution $W$ by estimating the norm $\ml^2_{e^{(s-s_0)/2} \xi}(Z(\cdot,s))$ by the norm $\av Z(\cdot,s_0) \av_{L^2_{\rho}}$ for $s_0$ large enough. This is done in Proposition \ref{pro:final} below. In this proof we will need the constructed semigroup estimates from the previous section.

Let us first formulate some auxiliary results. The next Proposition is merely a simple change of variables but it demonstrates how we are able to move from the $\ml^2_{\xi}$ norm to the $\mn^2_{|\xi|}$ norm.

\begin{pro} \label{pro:wpro1}
Let $Z$ be as above. If $s_0>1$ and $s' > 0$ and
$$
\sup_{s \in (s_0,s_0+s')} \ml^2_{e^{(s-s_0)/2} \xi}(Z(\cdot,s+\tau)) \le C_1,
$$
for every $|\xi| = 1$ and for every $\tau \ge 0$, then 
$$
\sup_{s \in (s_0,s_0+s')} \sup_{|\xi| \le 1} \ml^2_{e^{(s-s_0)/2}\xi}(Z(\cdot,s)) < C(C_1).
$$
\end{pro}
\emph{Proof.} For any $s \in (s_0,s_0+s')$, let $\xi(s) \in \mbr^N$ be such that $|\xi(s)| \le 1$ and
$$
\sup_{|\xi|\le 1} \ml^2_{e^{(s-s_0)/2}\xi}(Z(\cdot,s)) = \ml^2_{e^{(s-s_0)/2} \xi(s)}(Z(\cdot,s)).
$$
For $s \in (s_0,s_0+s')$ define a function $\beta$ through $e^{(\beta(s) -s_0)/2} = e^{(s-s_0)/2}|\xi(s)|$. This gives that $\beta(s) = s + 2\log(|\xi(s)|) \le s$. Let $I = \{s \in(s_0,s_0+s') \, : \, \beta(s) > s_0\}$.

Then, for $s \in I$ and for $\wh{\xi}(s) = \xi(s)/|\xi(s)|$, we have that
\begin{multline*}
\sup_{|\xi|\le 1} \ml^2_{e^{(s-s_0)/2}\xi}(Z(\cdot,s)) = \ml^2_{e^{(\beta(s)-s_0)/2}\wh{\xi}(s)}(Z(\cdot,s))
\\
= \ml^2_{e^{(\beta(s)-s_0)/2}\wh{\xi}(s)}(Z(\cdot,\beta(s) + s - \beta(s))) \le C_1,
\end{multline*}
by assumption since $|\wh{\xi}(s)| = 1$ and $s-\beta(s) \ge 0$.

Consider then $s \in (s_0,s_0+s') \setminus I$. Because $\beta(s) \le s_0$, we can write
\begin{multline*}
\sup_{|\xi| \le 1} \ml^2_{e^{(s-s_0)/2} \xi}(Z(\cdot,s)) = \ml^2_{e^{(s-s_0)/2} \xi(s)} (Z(\cdot,s)) 
\\
\le \mn^2_{e^{(\beta(s)-s_0)/2}}(Z(\cdot,s)) \le C_2 \left(\av Z(\cdot,s-1) \av_{L^2_{\rho}} + e^{-s} \right) \le C_3,
\end{multline*}
by using (\ref{esti3}) and (\ref{mn_Z_esti}) which finishes the proof. \nelio

In the next proposition we consider another type of regularizing property of the semigroup generated by the operator $A$. It is an $L^2$ - $L^{\infty}$ regularization for solutions with bounded gradient.

\begin{pro} \label{pro:boundedgrad2}
Let $Z$ be as above. Then
$$
Z(\xi + \gamma,s) \le C\Big( \ml^2_{e^{-1/2}\xi}(Z(\cdot,s-1))+ |\gamma|+e^{-s}\Big),
$$
for any $\gamma$ in $\mbr^N$, $|\xi| \le e^{s/2}R_1/2$ and $s \ge 1$.
\end{pro}

\emph{Proof.} Using the variation of constants formula together with the inequality (\ref{Z_equation_2}) and (\ref{ineq:Achi_1}), we get that, for every $|\xi| \le e^{s/2}R_1/2$,
\begin{multline*}
Z(\xi,s)
\\
\le e^{(A+\Gamma+A_1)\cdot 1}Z(\xi,s-1) + A_0\int_{s-1}^s \sup_{|\xi| \le e^{s/2}R_1/2} (e^{(A+\Gamma+A_1)(s-t)} \chi(\xi,t)\ud t
\\
\le e^{\Gamma+A_1} e^{A \cdot 1}Z(\xi,s-1) +A_0 C_1 e^{(\Gamma+A_1)s} \int_{s-1}^s e^{-(\Gamma+A_1)t} e^{-t}  \ud t.
\end{multline*}
By the representation formula (\ref{eq:semigroup1}) we estimate
\begin{multline*}
(e^{A\cdot 1}Z)(\xi,s-1) = \frac{1}{(1-e^{-1})^{N/2}}\int_{\mbr^N} \exp\left(-\frac{(e^{-1/2}\xi - \lambda)^2}{4(1-e^{-1})}\right) |Z(\lambda,s-1)| \ud \lambda
\\
\le \frac{1}{(1-e^{-1})^{N/2}} \Big\{ \int_{\mbr^N} e^{-\frac{|\eta|^2}{2(1-e^{-1})}} e^{\frac{|\eta|^2}{4}} \ud \eta \Big\}^{1/2} \cdot \Big\{ \int_{\mbr^N} |Z(e^{-1/2}\xi + \eta, s-1)|^2 e^{-\frac{|\eta|^2}{4}} \ud \eta \Big\}^{1/2}
\\
= C_2 \ml^2_{e^{-1/2}\xi}(Z(\cdot,s-1)).
\end{multline*}
Therefore, by using (\ref{nablaZ}), we have for every $|\xi| \le e^{s/2}R_1/2$ that
$$
Z(\xi+\gamma,s) \le Z(\xi,s) + A_4 |\gamma| \le  C_3 \left( \ml^2_{e^{-1/2}\xi}(Z(\cdot,s-1)) + |\gamma| + e^{-s} \right) 
$$
which gives us the claim. \nelio

In the next Proposition, we estimate the shifted $L^{\beta}$ -norm of $h(t)$ for some $\beta > 1$. What we want to obtain is that the norm is integrable with respect to $t$ if the shifted $L^2$-norm of $Z$ is bounded.

\begin{pro} \label{pro:hpro}
Let $Z$ be as above and let $R_1 e^{s_0/2} > 2$ and $s_0 >1$. If
$$
\sup_{s \in (s_0,s_0+s')} \ml^2_{e^{(s-s_0)/2} \xi}(Z(\cdot,s+\tau)) \le B,
$$
for every $\xi$ such that $|\xi| = 1$ and for every $\tau \ge 0$, then for any $\beta > 1$ there exist constants $M$, $\theta$, $C$ and $C'$ such that
\begin{multline*}
\ml^{\beta}_{e^{(t-s_0)/2} \xi} (h(\cdot,t)) \\
\le C \left( e^{-(t-s_0)/4\beta + C'\sqrt{t-s_0}} \ml^2_{e^{(t-s_0-\theta)/2}\xi}(Z(\cdot,t-\theta))^2 \right.
\\
\left.+  e^{-t}\ml^2_{e^{(t-s_0-\theta)/2}\xi}(Z(\cdot,t-\theta)) + e^{-t} \right),
\end{multline*}
for every $|\xi| = 1$ and $t \in (s_0+M,s_0+s')$, provided that $s'>M$. Here $M$ and $\theta$ depend only on the constants $\beta$, $\Gamma$, $R_1$, $N$ and $\{A_i\}_{i=0}^4$. The constants $C$ and $C'$ may, however, also depend on $B$.
\end{pro}

\emph{Proof.} By Proposition \ref{pro:wpro1} we have that
$$
\sup_{s \in (s_0,s_0+s')} \sup_{|\xi| \le 1} \ml^2_{e^{(s-s_0)/2}\xi}(Z(\cdot,s)) \le C_1(B).
$$
Since our assumptions imply that $|e^{(s-s_0)/2}\xi| < e^{s/2}R_1/2$ for every $|\xi| \le 1$, Proposition \ref{pro:boundedgrad2} tells us that
\begin{multline} \label{eq:hpro2}
\sup_{s \in (s_0+1,s_0+s')} \sup_{|\xi| \le 1} |Z(e^{(s-s_0)/2}\xi + \gamma,s)| 
\\
\le C_2'\sup_{s \in (s_0+1,s_0+s')} \sup_{|\xi| \le 1} \left(\ml^2_{e^{(s-s_0-1)/2}\xi}(Z(\cdot,s-1)) + |\gamma| + e^{-s}\right)
\\
< C_2(B) (1+|\gamma| ).
\end{multline}

Assuming that $M$ is large enough such that $R_1^{-1} e^{-M/2}M^{1/2} < \frac{1}{2}$, we have that
$$
e^{(t-s_0)/2} + (t-s_0)^{1/2} = R_1 e^{t/2} ( R_1^{-1} e^{-s_0/2} + R_1^{-1} e^{-t/2} (t-s_0)^{1/2})  < R_1 e^{t/2},
$$
for every $t> s_0 + M$. For $|y| > R_1 e^{t/2}$ and $|\xi| = 1$, we also have that $|y-e^{(t-s_0)/2}\xi| > R_1 e^{t/2}/2$, since we are assuming that $e^{-s_0/2} < R_1/2$.

Then, using the assumptions (\ref{esti1}) and (\ref{esti2}), we can estimate, for every $|\xi| = 1$ and $1 < M < t-s_0 < s'$, to obtain
\begin{multline*}
\ml^{\beta}_{e^{(t-s_0)/2}\xi}(h(\cdot,t))^{\beta}
\\
\le 2^{-\beta} \int_{|y| < e^{(t-s_0)/2}+ (t-s_0)^{1/2}} e^{\beta (\varphi(y) + \vartheta(y,t))} Z(y,t)^{2\beta} e^{-|y-e^{(t-s_0)/2}\xi|^2/4} \ud y +
\\
+ A_2^{\beta} \int_{e^{(t-s_0)/2}+(t-s_0)^{1/2}<|y|<R_1e^{t/2}} Z(y,t)^{2\beta} e^{-|y-e^{(t-s_0)/2}\xi|^2/4} \ud y
\\
+ A_1^{\beta} \int_{|y| > R_1 e^{t/2}} Z(y,t)^{\beta} e^{-|y-e^{-(t-s_0)/2}\xi|^2/4} \ud y
\end{multline*}
Above $\vartheta(y,s) \in (0,\wt{W}(y,s))$. Thus, by using (\ref{esti2}) and (\ref{eq:hpro2}) with $\gamma = (t-s_0)^{1/2} > M^{1/2} > 1$, we get that
$$
\vartheta(y,t) \le \max\{0,\wt{W}(y,t)\} \le 2C_2(B)(t-s_0)^{1/2},
$$
for every $|y| < e^{(t-s_0)/2}+ (t-s_0)^{1/2}$ and $t \in (s_0+M,s_0+s')$. This together with Hölder's inequality gives us the following estimate for the above expression
\begin{multline*}
\le 2^{-\beta} e^{2 \beta C_2(B) (t-s_0)^{1/2}} \ml^2_{e^{(t-s_0)/2}\xi}(e^{\beta \varphi}) \ml^2_{e^{(t-s_0)/2}\xi}(Z(\cdot,t)^{2\beta}) 
\\
+ A_2^{\beta} \Big\{\int_{|y| > (t-s_0)^{1/2}} e^{-|y|^2/4} \ud y \Big\}^{1/2} \ml^2_{e^{(t-s_0)/2}\xi} (Z(\cdot,t)^{2\beta})
\\
+ A_1^{\beta} \Big\{\int_{|y| > R_1 e^{t/2}/2} e^{-|y|^2/4} \ud y \Big\}^{1/2} \ml^2_{e^{(t-s_0)/2}\xi} (Z(\cdot,t)^{\beta}).
\end{multline*}
Finally, by (\ref{eq:phiestimate}) and (\ref{ml_Z_esti}), we can estimate the above by
\begin{multline*}
\le C_3 e^{2 \beta C_2(B) (t-s_0)^{1/2}} e^{-\beta(t-s_0)} \ml^{4\beta}_{e^{(t-s_0)/2}\xi}(Z(\cdot,t))^{2\beta}
\\
+ C_4 e^{-(t-s_0)/4} \ml^{4 \beta}_{e^{(t-s_0)/2}\xi}(Z(\cdot,t))^{2\beta} + C_5 e^{-\beta t} \ml^{2\beta}_{e^{(t-s_0)/2}\xi}(Z(\cdot,t))^{\beta}
\\
\le C_6 e^{2 \beta C_2(B) (t-s_0)^{1/2} -(t-s_0)/4} \left( \ml^2_{e^{(t-s_0-\theta)/2}\xi}(Z(\cdot,t-\theta)) + e^{-t} \right)^{2\beta}
\\
+ C_7 e^{-\beta t} \left(\ml^2_{e^{(t-s_0-\theta)/2}\xi}(Z(\cdot,t-\theta)) + e^{-t}\right)^{\beta},
\end{multline*}
for $\theta > 0$ satisfying $\frac{4\beta e^{-\theta}}{1+e^{-\theta}} < 1$. The claim follows after some simple estimations. \nelio

Now we are ready to state the proposition which is the cornerstone of the proof of Theorem \ref{theorem2}. At the end of this section, we will use Proposition \ref{pro:final} to obtain Theorem \ref{theorem2} as a relatively simple corollary.

\begin{pro} \label{pro:final}
Let $Z$ be as above. Then there exist constants $\ol{s}_0$, $C$, $K> 0$, depending only on $\{A_i\}_{i=0}^4$, $R_1$ and $\Gamma$, such that
$$
\ml^2_{e^{(s-s_0)/2} \xi}(Z(\cdot,s)) \le C \left( \sup_{s \ge s_0 - K} \av Z(\cdot,s) \av_{L^2_{\rho}} + e^{-s_0} \right),
$$
for every $|\xi| = 1$ and $s > s_0 \ge \ol{s}_0$.
\end{pro}

\emph{Proof.} Let $Z$ be as above, $|\xi|=1$ and define $Z_{\tau}(y,s) = Z(y,s+\tau)$ for $\tau \ge 0$. Then, since $\chi(s+\tau) \le \chi(s)$, we have that $Z_{\tau}$ satisfies the inequality (\ref{Z_equation}) with $h$ replaced by $h_{\tau}(y,s) = h(y,s+\tau)$. Therefore also (\ref{esti1})-(\ref{nablaZ}) hold for $Z_{\tau}$ and $h_{\tau}$ respectively, with the same constants $\{A_i\}_{i=1}^4$.

By the previous results, there exists a constant $M$ such that for every triplet $s$, $t$, $s_0$ for which $s > t > s_0  + M > M$, the following estimates (\ref{final_esti1})-(\ref{final_esti3}) hold. Firstly, we can assume that $M$ is large enough such that, by Proposition \ref{pro:semigroup2}, we have
\be \label{final_esti1}
\ml^2_{e^{(s-s_0)/2}\xi}(e^{\Lambda (s-s_0)}Z_{\tau}(\cdot,s_0)) \le C_1 \av Z_{\tau}(\cdot,s_0)\av_{L^2_{\rho}}, \qquad \text{for all } s \ge s_0 +M,
\ee
and, by Proposition \ref{pro:semigroup1},
\begin{multline}\label{final_esti2}
\ml^2_{e^{(s-s_0)/2}\xi}(e^{\Lambda (s-t)}h_{\tau}(\cdot,t)) \\
\le C_2(1-e^{-(s-t)})^{-\eps} \ml^{\beta}_{e^{(t-s_0)/2}\xi}(h_{\tau}(\cdot,t)), \quad \text{for } t \in (s_0 + M, s) \text{ and } s > s_0 +M,
\end{multline}
for some large $\beta$ and $\eps \in (0,1)$. By Corollary \ref{coro:semigroup1} and (\ref{esti1}),
\be \label{final_esti3}
\ml^2_{e^{(s-s_0)/2}\xi}(e^{\Lambda (s-t)}h_{\tau}(\cdot,t)|) \le C_3'(M) \av h_{\tau}(\cdot,t) \av_{L^2_{\rho}} \le C_3(M) \av Z_{\tau}(\cdot,t) \av_{L^2_{\rho}},
\ee
for $t \in [s_0,s_0+M]$ and $s > s_0 + K(M)$ when $K(M)>M$ is large enough.

Then, by (\ref{mn_Z_esti}) we have
\begin{multline} \label{final_esti4}
\ml^2_{e^{(s-s_0)/2}\xi}(Z_{\tau}(\cdot,s)) \le \mn^2_{e^{K(M)/2}}(Z_{\tau}(\cdot,s))
\\
\le C_4(M) \left( \av Z_{\tau}(\cdot, s-K(M)) \av_{L^2_{\rho}} + e^{-s} \right)
\end{multline}
for every $s_0 \ge K(M)$ and $s \in (s_0, s_0 + K(M))$.

Notice that in all these estimates the constants are independent of $\tau \ge 0$ and depend only on parameters such as $\{A_i\}_{i=0}^4$, $\Gamma$, $\beta$, $N$ and $M$. Since all these paratemers, apart from $M$, are fixed, we only point out the dependencies on $M$ below.

Let us estimate the $L^2_{e^{(s-s_0)/2}\xi}$ -norm of $Z(\cdot,s)$ for $s > s_0 + K(M) > s_0+M$ with $s_0 \ge K(M)$. The variation of constants formula and the inequality (\ref{Z_equation}) give that
\begin{multline*}
\ml^2_{e^{(s-s_0)/2}\xi}(Z_{\tau}(\cdot,s)) \le \ml^2_{e^{(s-s_0)/2}\xi}(e^{\Lambda (s-s_0)}Z_{\tau}(\cdot,s_0))
\\
+ \int_{s_0}^{s} \ml^2_{e^{(s-s_0)/2}\xi}(e^{\Lambda (s-t)}h_{\tau}(\cdot,t)) \ud t + \int_{s_0}^{s} \ml^2_{e^{(s-s_0)/2}\xi}(e^{\Lambda (s-t)}\chi(\cdot,t)) \ud t
\\
= T_1 + T_2 + T_3.
\end{multline*}
First we notice that (\ref{final_esti1}) gives an estimate for the term $T_1$. Then we split the integral in $T_2$ in two parts and use (\ref{final_esti2}) and (\ref{final_esti3}) to obtain
\begin{multline*}
T_2 = \int_{s_0}^{s_0+M} \ml^2_{e^{(s-s_0)/2}\xi}(e^{\Lambda (s-t)}h_{\tau}(\cdot,t)) \ud t + \int_{s_0+M}^s \ldots \ud t
\\
\le \int_{s_0}^{s_0+M}C_3(M) \av Z_{\tau}(\cdot,t) \av_{L^2_{\rho}} \ud t + \int_{s_0+M}^{s} C_2(1-e^{-(s-t)})^{-\eps} \ml^{\beta}_{e^{(t-s_0)/2}\xi}(h_{\tau}(\cdot,t)) \ud t
\\
\le C_5(M)  \sup_{s \ge s_0} \av Z_{\tau}(\cdot,s) \av_{L^2_{\rho}} + \int_{s_0+M}^{s} C_2(1-e^{-(s-t)})^{-\eps} \ml^{\beta}_{e^{(t-s_0)/2}\xi}(h_{\tau}(\cdot,t)) \ud t.
\end{multline*}
Similar estimates for the term $T_3$ imply
\begin{multline*}
T_3 \le \int_{s_0}^{s_0+M}C_3(M) \av \chi(\cdot,t) \av_{L^2_{\rho}} \ud t + \int_{s_0+M}^{s} C_2(1-e^{-(s-t)})^{-\eps} \ml^{\beta}_{e^{(t-s_0)/2}\xi}(\chi(\cdot,t)) \ud t
\\
\le C_6(M) \left( \int_{s_0}^{s_0+M} e^{-t} \ud t + \int_{s_0+M}^s (1-e^{-(s-t)})^{-\eps} e^{-t} \ud t \right) \le C_7(M) e^{-s_0},
\end{multline*}
which holds because we may assume that $|e^{(t-s_0)/2}\xi| \le e^{t/2}R_1/2$ and since we have that
$$
\ml^q_{\xi}(\chi(\cdot,t)) \le C e^{-t},
$$
for every $q \ge 1$ and $|\xi| \le e^{t/2}R_1/2$.

Therefore, we have proved that
\begin{multline} \label{eq:Z_variation}
\ml^2_{e^{(s-s_0)/2}\xi}(Z_{\tau}(\cdot,s))
\\
\le C_8(M) \left( \sup_{s \ge s_0} \av Z_{\tau}(\cdot,s) \av_{L^2_{\rho}} + e^{-s_0} \right) + \int_{s_0+M}^{s} C_2(1-e^{-(s-t)})^{-\eps} \ml^{\beta}_{e^{(t-s_0)/2}\xi}(h_{\tau}(t)) \ud t.
\end{multline}

Let $B = \max\{C_4(M),C_8(M)\} \ge 1$ and define, for every $s_0 \ge K(M)$,
\begin{multline*}
\ol{s}(s_0)
\\
 = \sup\Bigg\{ s' \ge s_0 : \ml^2_{e^{(s-s_0)/2}\xi}(Z_{\tau}(\cdot,s)) \le 2B \left( \sup_{s \ge s_0-K(M)} \av Z(\cdot,s) \av_{L^2_{\rho}} + e^{-s_0/2} \right),
\\
\text{ for every } s \in (s_0,s') \text{ and } \tau > 0 \Bigg\}.
\end{multline*}
Notice that because of the inequality (\ref{final_esti4}) we have $\ol{s}(s_0) \ge s_0+K(M)$. We want to show that there exists $\ol{s}_0$ such that $\ol{s}(s_0) = \infty$ whenever $s_0 \ge \ol{s}_0$.

By the previous Proposition, we may assume that $M$ is large enough such that
\begin{multline} \label{eq:hbeta}
\ml^{\beta}_{e^{(t-s_0)/2} \xi} (h_{\tau}(t)) \le C_9(M) \left( e^{-(t-s_0)/4\beta + C_9(M)\sqrt{t-s_0}} \ml^2_{e^{(t-s_0-\theta)/2}\xi}(Z_{\tau}(\cdot,t-\theta))^2 \right.
\\
\left. + e^{-t}\ml^2_{e^{(t-s_0-\theta)/2}\xi}(Z_{\tau}(\cdot,t-\theta)) + e^{-t} \right) ,
\end{multline}
for every $s_0+M < t < \ol{s}(s_0)$ and $e^{s_0/2} R_1 > 2$ and for some $\theta < M$.

We may also choose $K(M)$ to be such that all the previous inequalities in this proof hold and in addition $\frac{2 \beta e^{-K(M)}}{1 + e^{-K(M)}} < 1$. Then we can use Hölder's inequality and (\ref{ml_Z_esti}) to verify that, for $t \in (\ol{s}(s_0),\ol{s}(s_0) + \delta)$ and $\delta \in (0,1)$, we have
\begin{multline*}
\ml^{\beta}_{e^{(t-s_0)/2}\xi}(Z_{\tau}(t))^{\beta} = \int_{\mbr^N} |Z_{\tau}(t)|^{\beta} e^{-|y-e^{(t-s_0)/2}\xi|^2/4} \ud y
\\
\le \av Z_{\tau}(t)\av_{L^{2 \beta}_{\rho}}^{\beta} \left( \int_{\mbr^N} e^{y^2/4} e^{-|y-e^{(t-s_0)/2}\xi|^2/2} \ud y \right)^{1/2}
\\
\le C_{10}(M)^{\beta} \left( \av Z_{\tau}(\cdot,t-K(M))\av_{L^2_{\rho}} + e^{-t} \right)^{\beta} \left( \int_{\mbr^N} e^{-y^2/4} e^{|y|e^{(\ol{s}+1-s_0)/2}} \ud y \right)^{1/2}
\\
\le C_{10}(M)^{\beta}\left(\sup_{s \ge s_0} \av Z_{\tau}(\cdot,s) \av_{L^2_{\rho}} + e^{-s_0} \right)^{\beta} g(\ol{s}-s_0)^{\beta},
\end{multline*}
where the function $g$ is defined by the last equality and $\ol{s} = \ol{s}(s_0)$. Therefore,
\begin{multline} \label{eq:Zbeta_g}
\int_{\ol{s}}^{\ol{s} + \delta} (1-e^{-(\ol{s} + \delta-t)})^{-\eps} \ml^{\beta}_{e^{(t-s_0)/2}\xi}(Z_{\tau}(t)) \ud t 
\\
\le C_{10}(M) \delta^{1-\eps}g(\ol{s}-s_0) \left(\sup_{s \ge s_0}\av Z_{\tau}(\cdot,s) \av_{L^2_{\rho}} + e^{-s_0} \right).
\end{multline}


Now we use (\ref{eq:hbeta}) to estimate the integral term in (\ref{eq:Z_variation}) with $s = \ol{s} + \delta$, to get 
\begin{multline*}
I = \int_{s_0+M}^{\ol{s}+\delta} (1-e^{-(\ol{s} +\delta-t)})^{-\eps} \ml^{\beta}_{e^{(t-s_0)}\xi}(h_{\tau}(\cdot,t)) \ud t
\\
\le C_9 \int_{s_0+ M}^{\ol{s}} (1-e^{-(\ol{s} +\delta-t)})^{-\eps} e^{-(t-s_0)/4\beta + C_9\sqrt{t-s_0}} \ml^2_{e^{(t-\theta-s_0)/2}\xi}(Z_{\tau}(\cdot,t-\theta))^2 \ud t 
\\
+ C_{9} \int_{s_0+ M}^{\ol{s}} (1-e^{-(\ol{s} +\delta-t)})^{-\eps} e^{-t} \ml^2_{e^{(t-\theta-s_0)/2}\xi}(Z_{\tau}(\cdot,t-\theta)) \ud t 
\\
+ C_{9} \int_{s_0+ M}^{\ol{s}} (1-e^{-(\ol{s} +\delta-t)})^{-\eps} e^{-t} \ud t 
\\
+ C_{10} \delta^{1-\eps}g(\ol{s}-s_0) \left(\sup_{s \ge s_0} \av Z_{\tau}(\cdot,s) \av_{L^2_{\rho}} + e^{-s_0} \right).
\end{multline*}
Then use the definition of $\ol{s}$ in the two first integrals above to estimate
\begin{multline*}
I \le  2B \left( \sup_{s \ge s_0-K(M)} \av Z(\cdot,s) \av_{L^2_{\rho}} + e^{-s_0/2} \right)
\\
\cdot \Bigg\{2 C_9B \int _{s_0+ M}^{\ol{s}} (1-e^{-(\ol{s}+\delta-t)})^{-\eps} e^{-(t-s_0)/4\beta + C_9\sqrt{t-s_0}} \ud t \, \left( \sup_{s \ge s_0-K(M)} \av Z(\cdot,s) \av_{L^2_{\rho}} + e^{-s_0/2} \right)
\\
+C_{9} \int _{s_0+ M}^{s_0+\ol{s}} (1-e^{-(\ol{s}+\delta-t)})^{-\eps} e^{-t} \ud t  + C_{10}(M) \delta^{1-\eps}g(\ol{s}-s_0) \Bigg\}
\\
+ C_{9} \int_{s_0+ M}^{\ol{s}} (1-e^{-(\ol{s} +\delta-t)})^{-\eps} e^{-t} \ud t.
\end{multline*}
Defining the constant $C_{11}$ through
\begin{multline*}
\int_{s_0+M}^{\ol{s}} (1-e^{-(\ol{s}+\delta-t)})^{-\eps} e^{-(t-s_0)/4\beta + C_9\sqrt{t-s_0}} \ud t 
\\
\le C_{11}(M)'\int_{\ol{s}+\delta-1}^{\ol{s}} (1-e^{-\ol{s}+\delta-t})^{-\eps} \ud t + (1-e^{-1})^{-\eps}\int_{s_0+M}^{\ol{s}+\delta-1} e^{-(t-s_0)/4\beta + C_9\sqrt{t-s_0}} \ud t 
\\
\le  C_{11}(M),
\end{multline*}
and similarly
$$
\int_{s_0+ M}^{\ol{s}} (1-e^{-(\ol{s} +\delta-t)})^{-\eps} e^{-t} \ud t \le C_{12} e^{-s_0},
$$
we get that
\begin{multline} \label{I_esti1}
I \le 2B \left(\sup_{s \ge s_0-K(M)} \av Z(\cdot,s) \av_{L^2_{\rho}} +e^{-s_0/2} \right) 
\\
 \cdot\Bigg\{ 2 C_9 B C_{11}\left(\sup_{s \ge s_0-K(M)} \av Z(\cdot,s) \av_{L^2_{\rho}} + e^{-s_0/2} \right) + C_{9}C_{12} e^{-s_0} +  C_{10}(M) \delta^{1-\eps}g(\ol{s}-s_0)\Bigg\}
 \\
 + C_{9}C_{12}e^{-s_0}.
\end{multline}
By the assumption (\ref{esti3}) we may take $\ol{s}_0$ such that
\begin{equation*}
 2C_2 C_9 B  C_{11} \left( \sup_{s \ge s_0-K(M)} \av Z(\cdot,s) \av_{L^2_{\rho}} + e^{-s_0/2} \right) \le \frac{1}{16},
\end{equation*}
and
\be \label{s_0_esti_2}
C_2C_{9}C_{12}e^{-s_0/2} \le \frac{1}{16},
\ee
for every $s_0 > \ol{s}_0$. Then, by (\ref{eq:Z_variation}) and (\ref{I_esti1}), one has that 
\begin{multline} \label{ml_Z_final_ineq}
\ml^2_{e^{(\ol{s} +\delta)/2}\xi}(Z_{\tau}(\ol{s} +\delta)) \le 2B \left(\sup_{s \ge s_0-K(M)} \av Z(\cdot,s) \av_{L^2_{\rho}} +e^{-s_0/2} \right) 
\\
\cdot\Bigg\{ \frac{1}{2} + 2 C_2C_9 B C_{11}\left(\sup_{s \ge s_0-K(M)} \av Z(\cdot,s) \av_{L^2_{\rho}} + e^{-s_0/2} \right) 
\\
 + C_2C_{9}C_{12} e^{-s_0} +  C_2C_{10} \delta^{1-\eps}g(\ol{s}-s_0)\Bigg\}
 \\
 + C_2C_{9}C_{12}e^{-s_0}.
\\
\le 2B \left(\sup_{s \ge s_0-K(M)} \av Z(\cdot,s) \av_{L^2_{\rho}} +e^{-s_0/2} \right) 
\\
\left\{ \frac{1}{2} + \frac{1}{16} + \frac{1}{16} +C_2C_{10} \delta^{1-\eps}g(\ol{s}-s_0) + \frac{C_2C_{9}C_{12}}{2B} e^{-s_0/2} \right\}.
\end{multline}

Assuming that $\ol{s}(s_0) < \infty$ for some $s_0 > \ol{s}_0$ and defining $\delta= \delta(\ol{s}(s_0))$ to be small enough such that
$$
C_2C_{10}\delta^{1-\eps}g(\ol{s}-s_0) < \frac{1}{16},
$$
inequalities (\ref{s_0_esti_2}) and (\ref{ml_Z_final_ineq}) yield that
$$
\ml^2_{e^{(\ol{s} +\delta)/2}\xi}(Z_{\tau}(\ol{s} +\delta)) \le 2B \left( \sup_{s \ge s_0-K(M)} \av Z(\cdot,s) \av_{L^2_{\rho}} + e^{-s_0/2} \right)\ \cdot \frac{3}{4}
$$
for every $\tau > 0$. This is in contradiction with the definition of $\ol{s} = \ol{s}(s_0)$. We have thus proved the claim and reached the end of the proof of Proposition \ref{pro:final}. \nelio

\emph{Proof of Theorem \ref{theorem2}.}

Proposition \ref{pro:boundedgrad2} with $\gamma = 0$ gives us that, for $e^{-s_0/2}|\xi| \le \frac{R_1}{2}$,
$$
Z(e^{(s-s_0)/2}\xi,s) \le C \left( L^2_{e^{(s-s_0-1)/2}\xi}(Z(\cdot,s-1)) + e^{-s} \right),
$$
for every $s-1 >s_0$. Therefore, Proposition \ref{pro:final} states that
$$
Z(e^{(s-s_0)/2}\xi,s) \le C \left( \sup_{s \ge s_0 - K(M)} \av Z(\cdot,s) \av_{L^2_{\rho}} + e^{-s_0} \right),
$$
whenever $s_0>\ol{s}_0$, for some $\ol{s}_0$ large enough.

So, for $s$ large, $s_0>\ol{s}_0$ and $e^{-s_0/2}|\xi| \le R_1/2$, we have
\begin{multline*}
Z(e^{(s-s_0)/2}\xi,s)  = | -s + u(e^{-s_0/2}\xi, T-e^{-s}) - \varphi(e^{(s-s_0)/2}\xi) | \\
=  | -s + u(e^{-s_0/2}\xi, T-e^{-s}) +2\log(e^{(s-s_0)/2}\xi) -C_{\alpha} + o_{s \to \infty}(1) |  \\
= | u(e^{-s_0/2}\xi, T-e^{-s}) +2\log(e^{-s_0/2}\xi) -C_{\alpha} +  o_{s \to \infty}(1) |
\\
\le C \left( \sup_{s \ge s_0 - K(M)} \av Z(\cdot,s) \av_{L^2_{\rho}} + e^{-s_0} \right),
\end{multline*}
where $C_{\alpha}$ is as in Theorem \ref{theo:FP}.

By the above inequality, origin is the only blow-up point in the ball $B_{e^{-\ol{s}_0}}(0)$. Thereby, $u(x,t)$ is bounded in every set $B_{r_2}(0) \setminus B_{r_1}(0)$ for $0<r_1<r_2<e^{-\ol{s}_0}$ and parabolic estimates imply that also $u_t(x,t)$ is bounded in every such set. This gives us the existence of the limit $\lim_{t \to T} u(x,t)$ for $x \in B_{e^{-\ol{s}_0}}(0) \setminus\{0\}$.

By taking the limit as $s \to \infty$, we get
\bdm
| u(e^{-s_0/2}\xi, T) +2\log(e^{-s_0/2}\xi) -C_{\alpha} |\le C \left( \sup_{s \ge s_0 - K(M)} \av Z(\cdot,s) \av_{L^2_{\rho}} + e^{-s_0} \right).
\edm
Using the assumption (\ref{esti3}), we obtain
\bdm
\lim_{x \to 0} | u(x,T) + 2\log|x| - C_{\alpha} | = \lim_{s_0 \to \infty} | u(e^{-s_0/2}\xi, T) +2\log(e^{-s_0/2}\xi) -C_{\alpha} | = 0,
\edm
and thus we have found the blow-up profile. \nelio

\section{Revisiting the case of constant selfsimilar profile} \label{section_rev}

In this section we will briefly go through some results of the paper \cite{FP} and notice that the conclusions of Theorem \ref{theo:FP} above hold by assuming only that $u_0$ is radially symmetric and blow-up takes place at the origin with type I rate, thus verifying Theorem \ref{theorem_nondecrea} above.

In \cite{FP} we proved Theorem \ref{theo:FP} by showing that if $u$ is a radially symmetric and radially nonincreasing $L^1$-solution of equation (\ref{eq1}) with $f(u) = e^u$ on $[0,\mt]$ that blows up at $t = T < \mt$ and
$$
\log(T-t) + u(\sqrt{T-t}y,t) \to 0,
$$
uniformly on compact sets as $t \to T$, then either
\be \label{loglog_profile2}
\lim_{x \to 0} |u(x,T) +2\log|x| - \log|\log|x||| = C
\ee
or
\be \label{mlog_profile}
\lim_{x \to 0}|u(x,T) + m\log|x|| = C,
\ee
for some constant $C$. This will then imply that the blow-up is complete by Theorem 3.6 in \cite{Va}, thereby contradicting the assumption on $u$ being an $L^1$-solution on $[0,\mt]$. For an $L^1$-solution such as in Theorem \ref{theo:FP} the only possibility is thus a nonconstant selfsimilar blow-up profile, see details in \cite{FP}.

We will now demonstrate that we do not actually need to assume that $u$ is radially nonincreasing in order for this analysis to go through.

Let
$$
\wt{u}(x,t) = \zeta(x)u(x,t) - (\log(T-t) + 1)(1-\zeta(x))
$$
be a continuation of $u$ to the whole space $\mbr^N$, where $\zeta \in C^{\infty}(\mbr^N)$ and $\zeta(x) = 1$ for $|x| \le R_1< R_2 < R$ and $\zeta(x) = 0$ for $|x| > R_2$. It is proved in Propositions 3.4 and 3.6 in \cite{FP} that the following is true.
\begin{pro} \label{pro:FP_W_tau}
Let $u$ be a radially symmetric solution of (\ref{eq1}) that blows up with type I rate at $(x,t)=(0,T)$ and assume that the convergence (\ref{eq:ssconv1}) holds with $\varphi=0$. Then 
$$
\lim_{t \to T} \left( \log(T-t) + u(\lambda(t) y,t) \right) = -\log\left(1 + \sum_{|\alpha| = m}c_{\alpha} y^{\alpha} \right),
$$
uniformly for $y$ in compact sets, for some $m \ge 2$ and constants $c_{\alpha}$, where $\alpha$ is a multi-index and $|\alpha| = \alpha_1 + \ldots + \alpha_N$. The function $\lambda(t)$ is defined by $\lambda(t) = |\log(T-t)| \sqrt{T-t}$ if $m=2$, and by $\lambda(t) = (T-t)^{1/m}$ if $m > 2$.
\end{pro}
Notice that there is no assumption on $u$ being radially nonincreasing in the above Proposition.

We show that the above Proposition implies that either (\ref{loglog_profile2}) or (\ref{mlog_profile}) holds. Therefore the conclusion of Theorem \ref{theorem_nondecrea} is obtained by using Theorem 3.4 in \cite{Va} and an energy argument as in \cite{FP}. We will first define an auxiliary function $W_{\tau}$ and describe some of its properties below. Then we prove Proposition \ref{theorem_W_tau} below, which gives that the $L^2$-norm of $W_{\tau}(\cdot,s)$ is controlled by the $L^2$-norm of $W_{\tau}(\cdot,0)$. Profiles (\ref{loglog_profile2}) and (\ref{mlog_profile}) are obtained at the very end of this section.

For fixed $\xi$ with $|\xi| \le 1$, define
\begin{multline*}
W_{\tau}(y,s) = \log(T-\tau) + \wt{u}(\lambda(\tau) \xi + \sqrt{T-\tau}\sqrt{1-t} y,\tau + (T-\tau)t)
\\
+ \log \left(1-t +\sum_{|\alpha|=m} c_{\alpha} \xi^{\alpha} \right),
\end{multline*}
for $y \in \mbr^N$ and $s \ge 0$, where $\tau \in (0,T)$ and $s = -\log(1-t)$. Then, by using the above Proposition \ref{pro:FP_W_tau}, $W_{\tau}$ satisfies
\be \label{W_tau_to_0}
\av W_{\tau}(\cdot,0) \av_{L^2_{\rho}} \to 0 \quad \text{ and } \quad \av W_{\tau}(\cdot,0) \av_{L^{\beta}_{\rho}} \to 0,	
\ee
as $\tau \to T$ and, by (\ref{eq:nablau}), there exists a constant $A_0>0$ such that 
\be \label{W_tau_grad_esti}
|\nabla W_{\tau}(y,s)| \le A_0,
\ee
for every $(y,s) \in \mbr^N \times (0,\infty)$.  It can be also verified that $W_{\tau}$ solves
$$
(W_{\tau})_s = \Delta W_{\tau} -\frac{y}{2} \nabla W_{\tau} + e^{\wt{\phi}} W_{\tau} + f_{\tau} = A W_{\tau} + e^{\wt{\phi}} W_{\tau} + f_{\tau} ,
$$
where $A = \Delta - \frac{y}{2}\nabla$ and
$$
\wt{\phi}(s) = -s - \log\Bigg(e^{-s} + \sum_{|\alpha| = m} c_{\alpha} \xi^{\alpha}\Bigg).
$$
Above $f_{\tau}$ is a certain function, that we do not explicitly write out here, satisfying
\be \label{ftau_ineq}
|f_{\tau}| \le A_1 |W_{\tau}| \quad \text{and} \quad |f_{\tau}| \le A_2 |W_{\tau}|^2,
\ee
for some constants $A_1,A_2>0$ and
\be \label{f_tau_ineq1}
|f_{\tau}| \le e^{\wt{\phi}}(e^{W_{\tau}} - 1 - W_{\tau}),
\ee
for $|\lambda(\tau) \xi + \sqrt{T-\tau}\sqrt{1-t} y| \le R_1$. We have the following result.

\begin{pro} \label{theorem_W_tau}
Let $W_{\tau}$ be as above. There exist constants $C$ and $\ol{\tau}>0$ such that
$$
\av W_{\tau}(\cdot,s) \av_{L^2_{\rho}} \le C \av W_{\tau}(\cdot,0) \av_{L^2_{\rho}},
$$
for every $s \ge 0$ and $\tau > \ol{\tau}$.
\end{pro}

\emph{Proof.} Let $Z_{\tau} = |W_{\tau}|$. Since $e^{\wt{\phi}}$ is bounded and since $|f_{\tau}| \le A_1 Z_{\tau}$, we have a constant $C_1>0$ such that for any $s \ge s_0 \ge 0$ it holds
\be \label{Z_tau_esti0}
Z_{\tau}(\cdot,s) \le e^{(A+C_1) (s-s_0)} Z_{\tau}(\cdot,s_0).
\ee
By (\ref{W_tau_to_0}) one has that $\av Z_{\tau}(\cdot,0) \av_{L^2_{\rho}} \le C$, for some $C>0$ independent of $\tau$, which implies
\be \label{Z_tau_esti1}
\av Z_{\tau}(\cdot,s) \av_{L^2_{\rho}} \le e^{C_1 s_0} \av Z_{\tau}(\cdot,0) \av_{L^2_{\rho}} \le C_2(s_0),
\ee
for every $s \le s_0$ and $\tau \in (0,T)$ with $C_2(s_0)$ independent of $\tau$.

Now define
$$
\ol{s}(s_0) = \sup\{s'>0 \, : \, \av Z_{\tau}(\cdot,s)\av_{L^2_{\rho}} \le 4 e^{C_1s_0} \av Z_{\tau}(\cdot,0) \av_{L^2_{\rho}} \text{ for every } s \in (0,s')\}
$$
and notice that (\ref{Z_tau_esti1}) implies $\ol{s}(s_0) > s_0$.

Using the representation (\ref{eq:semigroup1}) and (\ref{Z_tau_esti0}) together with the definition of $\ol{s}(s_0)$, we obtain
\begin{multline} \label{Z_tau(0,s)_esti}
Z_{\tau}(0,s) \le e^{C_1} e^{A 1} Z_{\tau}(0,s-1)
\\
\le C_3' \int_{\mbr^N} e^{-\frac{|\lambda|^2}{4(1-e^{-1})}} |Z_{\tau}(\lambda,s-1)| \ud \lambda
\\
\le C_3' \left( \int_{\mbr^N} |Z_{\tau}(\lambda,s-1)|^2 e^{-|\lambda|^2/4} \ud y \right)^{1/2} \left(\int_{\mbr^N} e^{|\lambda|^2/4} e^{-|\lambda|^2/2} \ud \lambda \right)^{1/2}
\\
\le C_3 e^{C_1 s_0} \av Z_{\tau}(\cdot,0) \av_{L^2_{\rho}} \le C_4(s_0),
\end{multline}
for every $s \in (1,\ol{s}(s_0))$. For $|y| \le \sqrt{s}$ we have that $|\lambda(\tau) \xi + \sqrt{T-\tau} e^{-s/2} y| \le R_1$ for $\tau$ large enough, and for such $y$ the inequalities (\ref{W_tau_grad_esti}), (\ref{f_tau_ineq1}) and (\ref{Z_tau(0,s)_esti}) thus imply
\begin{multline} \label{f_tau_esti1}
|f_{\tau}(y,s)| \le \frac{1}{2} e^{\wt{\phi}(s) + Z_{\tau}(y,s)} Z_{\tau}(y,s)^2 
\\
\le C_5 e^{-s} e^{Z_{\tau}(0,s) + A_0\sqrt{s}} Z_{\tau}(y,s)^2 \le C_6(s_0) e^{-s + A_0 \sqrt{s}} Z_{\tau}(y,s)^2.
\end{multline}
Now we can estimate the $L^2$-norm of $f_{\tau}$, by using (\ref{ftau_ineq}), (\ref{f_tau_esti1}), Hölder's inequality, Proposition \ref{pro:velaz2} and the definition of $\ol{s}(s_0)$, to obtain
\begin{multline*}
\av f_{\tau}(\cdot,s) \av_{L^2_{\rho}}^2 \le C_6(s_0)^2 \int_{|y| \le \sqrt{s}}  e^{-2s + 2A_0 \sqrt{s}} Z_{\tau}(y,s)^4 e^{-|y|^2/4} \ud y
\\
+ \int_{|y| > \sqrt{s}} A_2^2 Z_{\tau}(y,s)^4 e^{-|y|^2/4} \ud y
\\
\le  C_6(s_0)^2 e^{-2s + 2A_0 \sqrt{s}} \av Z_{\tau}(\cdot,s) \av_{L^4_{\rho}}^4
\\
+ A_2^2 \left( \int_{|y| > \sqrt{s}} e^{-|y|^2/4} \ud y \right)^{1/2} \av Z_{\tau}(\cdot,s) \av_{L^8_{\rho}}^4
\\
\le C_7(s_0) e^{-2s + 2A_0 \sqrt{s}} \av Z_{\tau}(\cdot,s-K) \av_{L^2_{\rho}}^4 + C_8 e^{-s/8} \av Z_{\tau}(\cdot,s-K) \av_{L^2_{\rho}}^4 
\\
\le C_9(s_0)^2 e^{-s/8 + 2A_0 \sqrt{s}} \av Z_{\tau}(\cdot,0)\av_{L^2_{\rho}}^4,
\end{multline*}
for any $s \in (K,\ol{s}(s_0))$ and $s_0>K$.

Now, let $C_{10}$ be such that $e^{\wt{\phi}(s)} \le C_{10}e^{-s}$ and take $s_0>K$ large enough to satisfy
$$
C_{10} e^{-s_0} < \frac{1}{4}
$$
and define $\ol{\tau}=\ol{\tau}(s_0)$ to be such that
$$
\frac{C_9(s_0)}{4 e^{C_1s_0}} \int_{s_0}^{\infty} e^{-t/16 + A_0 \sqrt{t}} \ud t \av Z_{\tau}(\cdot,0) \av_{L^2_{\tau}} < \frac{1}{4},
$$
for every $\tau > \ol{\tau}$.

Then by using the variation of constants formula together with the previous estimates, we obtain 
\begin{multline*}
\av Z_{\tau}(\cdot,s) \av_{L^2_{\rho}}
\\
\le \av e^{A(s-s_0)} Z_{\tau}(\cdot,s_0) \av_{L^2_{\rho}} + \int_{s_0}^s \av e^{A(s-t)} e^{\wt{\phi}(t)} Z_{\tau}(\cdot,t) \av_{L^2_{\rho}} \ud t + \int_{s_0}^s \av e^{A(s-t)} f_{\tau}(\cdot,t) \av_{L^2_{\rho}} \ud t 
\\
\le  \av Z_{\tau}(\cdot,s_0) \av_{L^2_{\rho}}  + C_{10}4 e^{C_1 s_0} \int_{s_0}^s e^{-t} \ud t \av Z_{\tau}(\cdot,0)\av_{L^2_{\rho}} 
\\
+ C_9(s_0) \int_{s_0}^s e^{-t/16 + A_0 \sqrt{t}} \ud t \av Z_{\tau}(\cdot,0) \av_{L^2_{\rho}}^2 
\\
\le 4 e^{C_1 s_0} \av Z_{\tau}(\cdot,0) \av_{L^2_{\rho}} \left( \frac{1}{4} + C_{10} e^{-s_0} +  \frac{C_9(s_0)}{4 e^{C_1 s_0}} \int_{s_0}^{\infty} e^{-t/16 + A_0 \sqrt{t}} \ud t \av Z_{\tau}(\cdot,0)\av_{L^2_{\rho}} \right) 
\\
\le 3 e^{C_1 s_0} \av Z_{\tau}(\cdot,0) \av_{L^2_{\rho}},
\end{multline*}
for every $s \in (s_0,\ol{s}(s_0))$ and $\tau > \ol{\tau}$. This proves that $\ol{s}(s_0) = \infty$ and the claim follows. \nelio

Using the previous Proposition and (\ref{Z_tau(0,s)_esti}) we have that
$$
|W_{\tau}(0,s)| \le C \av W_{\tau}(\cdot,s-1) \av_{L^2_{\rho}} \le C \av W_{\tau}(\cdot,0) \av_{L^2_{\rho}},
$$
for every $\tau > \ol{\tau}$ and $s >1$, when $\ol{\tau}$ is close enough to $T$. Rewriting this gives
$$
\Big|\log(T-\tau) + u(\lambda(\tau) \xi, \tau + (T-\tau)(1-e^{-s})) + \log\Big(e^{-s} + \sum_{|\alpha|=m} c_{\alpha} \xi^{\alpha}\Big)\Big| \le C\av W_{\tau}(\cdot,0) \av_{L^2_{\rho}} ,
$$
and by taking the limit as $s \to \infty$, one notices that
$$
\Big|\log(T-\tau) + u(\lambda(\tau) \xi,T) + \log\Big(\sum_{|\alpha|=m} c_{\alpha} \xi^{\alpha}\Big) \Big| \le C\av W_{\tau}(\cdot,0) \av_{L^2_{\rho}} \to 0,
$$
as $\tau \to T$. The desired blow-up profiles (\ref{loglog_profile2}) and (\ref{mlog_profile}) are then obtained by a change of variables as in \cite{FP}.


\end{document}